\documentclass[a4paper,9pt]{article}
\usepackage{makeidx}
\usepackage{amsfonts}
\usepackage{vmargin}
\usepackage{amssymb}
\usepackage{amsmath}
\usepackage{graphicx}
\usepackage[latin1]{inputenc}
\usepackage{amsthm}
\usepackage{color}
\usepackage{mathrsfs}
\usepackage{fancyhdr}
\usepackage{setspace}
\usepackage{float}
\usepackage[center,small,bf,ruled]{caption}

\setmarginsrb{40mm}{27mm}{40mm}{27mm}
             {0mm}{10mm}{0mm}{10mm}

\numberwithin{equation}{section}

\begin{document}

\title{Fractional relaxation equations and Brownian crossing\\
probabilities of a random boundary}
\author{L.Beghin\thanks{%
Dep. of Statistical Sciences, Sapienza University of Rome, p.le A.Moro 5,
00185, Rome, Italy. E-mail address: luisa.beghin@uniroma1.it}}
\maketitle

\begin{abstract}
We analyze here different forms of fractional relaxation equations of order $%
\nu \in \left( 0,1\right) $ and we derive their solutions both in analytical
and in probabilistic forms. In particular we show that these solutions can
be expressed as crossing probabilities of random boundaries by various types
of stochastic processes, which are all related to the Brownian motion $B$.
In the special case $\nu =1/2$, the fractional relaxation is proved to
coincide with $\Pr \left\{ \sup_{0\leq s\leq t}B(s)<U\right\} $, for an
exponential boundary $U.$ When we generalize the distributions of the random
boundary, passing from the exponential to the Gamma density, we obtain more
and more complicated fractional equations.


\end{abstract}

\section{Introduction}

The following differential equation%
\begin{equation}
\frac{d}{dt}p(t)=-\lambda p(t),\quad t>0  \label{uno}
\end{equation}%
is known in the physics literature as the \textit{relaxation equation}. The
solution to (\ref{uno}), with initial condition $p(0)=1,$ is clearly equal
to $p(t)=e^{-\lambda t}.$ Since the end of the Nineties an intensive
research activity has been developed, aimed at the application of fractional
calculus to mathematical physics: many classical equations have been
modified by substituting the integer-order derivatives with the fractional
ones. Equation (\ref{uno}) has been extended in the following fractional
sense:%
\begin{equation}
\frac{d^{\nu }}{dt^{\nu }}\psi (t)=-\lambda \psi (t),\quad t>0  \label{due}
\end{equation}%
where $\nu \in \left( 0,1\right) $ and $\frac{d^{\nu }}{dt^{\nu }}$
represents the fractional derivative according to the Caputo definition, i.e.%
\begin{equation}
\frac{d^{\nu }}{dt^{\nu }}u(t)=\left\{
\begin{array}{l}
\frac{1}{\Gamma (m-\nu )}\int_{0}^{t}\frac{1}{(t-s)^{1+\nu -m}}\frac{d^{m}}{%
ds^{m}}u(s)ds\text{,\qquad for }m-1<\nu <m \\
\frac{d^{m}}{dt^{m}}u(t)\text{,\qquad for }\nu =m,%
\end{array}%
\right. ,  \label{tre}
\end{equation}%
with $m=\left\lfloor \alpha \right\rfloor +1.$ Obviously, for $\nu =1$ the
\textit{fractional relaxation equation} (\ref{due}) coincides with the
standard equation (\ref{uno}).

Equation (\ref{due}) has been studied in some papers, such as \cite{mai1},
\cite{mai} and its solution was given analytically in terms of the
Mittag-Leffler function as:%
\begin{equation}
\psi _{\nu }(t)=E_{\nu ,1}(-\lambda t^{\nu }),  \label{qua}
\end{equation}%
where
\begin{equation}
E_{\alpha ,\beta }(z)=\sum_{r=0}^{\infty }\frac{\,z^{r}}{\Gamma (\alpha
r+\beta )},\quad \alpha ,\beta \in \mathbb{C},\text{ }Re(\alpha ),%
Re(\beta )>0.  \label{ML}
\end{equation}%
The analysis of the fractional relaxation equation has mainly physical
motivations, for instance to study the electromagnetic properties of a wide
range of materials (which display a long memory, instead of exponential,
decay, see \cite{sib} and \cite{siba}) as well as the rheological models for
the description of some viscoelastic materials (see \cite{met}, \cite{glo},
\cite{non} and \cite{sci}).

Moreover, the so-called Mittag-Leffler distribution has been often applied
to statistics (for example in \cite{lin} and \cite{pil1}) or to queuing
theory in \cite{pil2}.

Actually the solution $\psi _{\nu }(t),t>0$ can be expressed in
probabilistic terms in two interesting forms, that we will present and
explore here. The first form represents the probability of no events up to
time $t$ (or survival probability), for the so-called \emph{fractional
Poisson process} $\mathcal{N}_{\nu }(t),t>0$ (see, among the others, \cite%
{las}, \cite{wan}, \cite{mai2}, \cite{beg1}, and \cite{beg2}). Indeed the
following equality holds
\begin{equation}
\psi _{\nu }(t)=p_{0}^{\nu }(t)=\Pr \left\{ \mathcal{N}_{\nu }(t)=0\right\}
\label{sei}
\end{equation}%
and thus we can apply to $\psi _{\nu }(t)$ the results obtained in the above
cited articles.\textit{\ }For example we will resort to the equality of the
one-dimensional distribution between $\mathcal{N}_{\nu }$ and a composition
of the standard Poisson process $N(t)$ with a random time-process $\mathcal{T%
}_{\nu }(t),$ i.e. $N(\mathcal{T}_{\nu }(t)),t>0.$ Thus, thanks to (\ref{sei}%
), we can write%
\begin{equation}
\psi _{\nu }(t)=\int_{0}^{\infty }e^{-\lambda y}q_{\nu }(y,t)dy=\Pr \left\{
\mathcal{T}_{\nu }(t)<U\right\} ,  \label{set}
\end{equation}%
where $q_{\nu }(y,t)$ is the density of $\mathcal{T}_{\nu }$ (which is
itself solution to a fractional diffusion equation) and $U$ is an
exponential random variable with parameter $\lambda >0$. Formula (\ref{set})
is particularly interesting in the special case where $\nu =1/2,$ since it
becomes%
\begin{equation}
\psi _{1/2}(t)=\int_{0}^{\infty }e^{-\lambda y}\frac{e^{-y^{2}/4t}}{\sqrt{%
\pi t}}dy=\Pr \left\{ |B(t)|<U\right\} ,  \label{ott}
\end{equation}%
where $B$ is a Brownian motion starting from zero and with variance $2t.$

As a consequence, a second probabilistic interpretation of the solution to
the fractional relaxation equation (\ref{due}) can be given in terms of
\textit{crossing probability} of a random boundary by a standard Brownian
motion, for $\nu =1/2.$ Indeed it is well known that the following
relationship holds:%
\begin{equation*}
\Pr \left\{ |B(t)|<z\right\} =\Pr \left\{ \sup_{0\leq s\leq t}B(s)<z\right\}
=\Pr \left\{ B(s)<z,\text{ }\forall \text{ }s\in \left( 0,t\right) \right\} ,
\end{equation*}%
where the last expression is commonly referred to as \emph{crossing
probability}.

For other values of $\nu ,$ an analogue result holds true, but for less
known processes, such as the iterated Brownian motion (for $\nu =1/2^{n}$)
or the \emph{Airy process}\ (for $\nu =1/3$).

Moreover, the expression (\ref{set}) shows that the solution to (\ref{due})
can be expressed as a standard relaxation with random time represented by $%
\mathcal{T}_{\nu }$, i.e. as $\psi (\mathcal{T}_{\nu }(t))$. The results
given in \cite{mer1} permit also to express the solution as a time-changed
relaxation via an inverse stable subordinator $E(t),$ i.e. as $\psi _{\nu
}(t)=\psi (E(t)).$ In fact $\psi (\mathcal{T}_{\nu }(t))$ and $\psi (E(t))$
share the one-dimensional distributions and therefore the two approaches can
be considered equivalent.

In the successive sections we analyze some extensions of the result (\ref%
{set}) in the following directions:

\begin{itemize}
\item We consider other random time-processes in place of $\mathcal{T}_{\nu
} $ and therefore in (\ref{ott}) instead of Brownian motion: for example,
the sojourn time of a Brownian motion on the positive half-line, the
first-passage time of a Brownian motion through a certain level, the elastic
Brownian motion (by analogy with the analysis carried out for the fractional
Poisson process in \cite{beg4}).

\item We consider a different random variable (i.e. the Gamma) instead of $U$
in (\ref{ott});

\item We introduce in (\ref{due}) an assumption of \emph{distributed
fractional derivative} (see \cite{mai}, \cite{beg3}).
\end{itemize}

\section{Fractional relaxation equation of order $\protect\nu $}

A first probabilistic expression of the solution $\psi _{\nu }(t),t>0$\ to
equation (\ref{due}) can be found by considering that the latter coincides
with the fractional equation satisfied by the survival probability (i.e. the
probability of no events up to time $t$) of a fractional Poisson process of
order $\nu \in \left( 0,1\right) .$ Let $\mathcal{N}_{\nu }(t),t>0$, denote
the process with probabilities $p_{k}^{\nu }(t)$ solving the following
recursive differential equation
\begin{equation}
\frac{d^{\nu }p_{k}^{\nu }}{dt^{\nu }}=-\lambda (p_{k}^{\nu }-p_{k-1}^{\nu
}),\quad k\geq 0,t>0,  \label{ai}
\end{equation}%
with initial conditions
\begin{equation}
p_{k}^{\nu }(0)=\left\{
\begin{array}{c}
1\qquad k=0 \\
0\qquad k\geq 1%
\end{array}%
\right.  \label{ai2}
\end{equation}%
and $p_{-1}^{\nu }(t)=0$. The process $\mathcal{N}_{\nu }$ has been studied
in a series of papers (for example in \cite{las}, \cite{mai2} and \cite{beg1}%
, in the homogeneous case, and in \cite{wan}, in the non-homogeneous case)
and its distribution has been expressed in analytic forms in terms of
derivatives of Mittag-Leffler function or as generalized Mittag-Leffler
(GML) functions
\begin{equation}
E_{\alpha ,\beta }^{\gamma }(z)=\sum_{j=0}^{\infty }\frac{\left( \gamma
\right) _{j}\,z^{j}}{j!\Gamma (\alpha j+\beta )},\quad \alpha ,\beta ,\gamma
\in \mathbb{C},\text{ }Re(\alpha ),Re(\beta ),Re(\gamma )>0,  \label{gml2}
\end{equation}%
where $\left( \gamma \right) _{j}=\gamma (\gamma +1)...(\gamma +j-1)$ (for $%
j=1,2,...,$ and $\gamma \neq 0$) and $\left( \gamma \right) _{0}=1$ (see
\cite{beg2}). Moreover in \cite{beg1} a probabilistic expression of the
process has been given, as composition of a\ standard Poisson process $N$
with a random time argument $\mathcal{T}_{\nu }$, independent of $N$. The
following equality in distribution was proved to hold:%
\begin{equation}
\mathcal{N}_{\nu }(t)\overset{i.d.}{=}N(\mathcal{T}_{\nu }(t)),  \label{ult}
\end{equation}%
where $\mathcal{T}_{\nu }(t)$ possesses transition density $q_{\nu }(y,t)$
coinciding with the folded solution to the fractional diffusion equation%
\begin{equation}
\frac{\partial ^{2\nu }v}{\partial t^{2\nu }}=\frac{\partial ^{2}v}{\partial
y^{2}},\qquad t>0,\text{ }y\in \mathbb{R},\quad v(y,0)=\delta
(y),v_{t}(y,0)=0  \label{fold}
\end{equation}%
i.e. with%
\begin{equation}
q_{\nu }(y,t)=\left\{
\begin{array}{l}
2v(y,t),\quad y\geq 0 \\
0,\quad y<0%
\end{array}%
\right. .  \label{qu}
\end{equation}%
Alternatively, it has been also proved in \cite{FP} and in \cite{mer1} that $%
q_{\nu }(y,t)$ solves the following equation%
\begin{equation}
\frac{\partial ^{\nu }q}{\partial t^{\nu }}=-\frac{\partial q}{\partial y}%
,\qquad t>0,\quad q(y,0)=\delta (y),  \label{nofold}
\end{equation}%
where, in this case, $y>0.$ In any case we can write
\begin{equation*}
p_{k}^{\nu }(t)=\Pr \left\{ \mathcal{N}_{\nu }(t)=k\right\} =\frac{(\lambda
t)^{k}}{k!}\int_{0}^{+\infty }e^{-\lambda y}q_{\nu }(y,t)dy,
\end{equation*}%
so that we immediately have, in view of (\ref{ai}) for $k=0,$ that%
\begin{equation}
\psi _{\nu }(t)=p_{0}^{\nu }(t)=\Pr \left\{ \mathcal{N}_{\nu }(t)=0\right\}
=\int_{0}^{+\infty }e^{-\lambda y}q_{\nu }(y,t)dy.  \label{ai3}
\end{equation}%
Therefore, in view of (\ref{ult}), the fractional relaxation $\psi _{\nu }$
can be expressed as composition of the standard relaxation with the random
time $\mathcal{T}_{\nu }$:%
\begin{equation*}
\psi _{\nu }(t)=\psi (\mathcal{T}_{\nu }(t)),\qquad t>0.
\end{equation*}

\subsection{Exponential boundary crossing probabilities of Brownian motion}

As a consequence of (\ref{ai3}) a second probabilistic form of the solution
in terms of boundary crossing probabilities is obtained in the following
result.

\

\noindent \textbf{Theorem 2.1 }\textit{Let }$U$\textit{\ be a random
boundary exponentially distributed (with parameter }$\lambda >0)$\textit{,
then the crossing probability of }$U$\textit{\ by the independent random
process }$\mathcal{T}_{\nu }(t)$\textit{\ with transition density }$q_{\nu
}(y,t),$\textit{\ i.e.}%
\begin{equation}
\psi _{\nu }(t)=\Pr \left\{ \mathcal{T}_{\nu }(t)<U\right\} ,  \label{pr}
\end{equation}%
\textit{satisfies the fractional relaxation equation (\ref{due}), with
initial condition }$\psi _{\nu }(0)=1.$

\noindent \textbf{Proof \ }We consider now the analytic expression of the
folded solution $q_{\nu }(y,t)$ to problem (\ref{fold}), in terms of the
Wright function%
\begin{equation*}
\mathcal{W}_{\alpha ,\beta }(x)=\sum_{j=0}^{\infty }\frac{x^{j}}{j!\Gamma
(\alpha j+\beta )},\qquad \alpha \geq -1,\;\beta >0,\;x\in \mathbb{R},
\end{equation*}%
which reads%
\begin{equation*}
q_{\nu }(y,t)=2v(y,t)=\frac{1}{t^{\nu }}W_{-\nu ,1-\nu }\left( -\frac{y}{%
t^{\nu }}\right) ,\quad y,t>0
\end{equation*}%
(see, for example, \cite{mai1}). Therefore we can rewrite (\ref{pr}) as
follows%
\begin{eqnarray}
\psi _{\nu }(t) &=&\Pr \left\{ \mathcal{T}_{\nu }(t)<U\right\}  \label{pr2}
\\
&=&\lambda \int_{0}^{\infty }e^{-\lambda y}\Pr \left\{ \mathcal{T}_{\nu
}(t)<y\right\} dy  \notag \\
&=&\frac{\lambda }{t^{\nu }}\int_{0}^{\infty }e^{-\lambda
y}\int_{0}^{y}W_{-\nu ,1-\nu }\left( -\frac{z}{t^{\nu }}\right) dzdy  \notag
\\
&=&\frac{1}{t^{\nu }}\int_{0}^{\infty }e^{-\lambda z}W_{-\nu ,1-\nu }\left( -%
\frac{z}{t^{\nu }}\right) dz  \notag \\
&=&E_{\nu ,1}(-\lambda t^{\nu }),  \notag
\end{eqnarray}%
by the well-known formula of the Laplace transform of the Wright function
(see \cite{po}, formula (1.165), p.39). The last expression in (\ref{pr2})
coincides with the solution to equation (\ref{due}) given in (\ref{qua}%
).\hfill $\blacksquare $

\

The previous results can be particularly relevant in the special case where $%
\nu =1/2$, since the random process $\mathcal{T}_{\nu }$ reduces to a
reflecting Brownian motion: indeed in this case the equation (\ref{fold})
governing the process coincides with the heat equation and $q_{1/2}(y,t)$
becomes the Gaussian with variance $2t$, folded with respect to the origin.
Therefore the fractional relaxation equation of order $1/2$ is solved by%
\begin{eqnarray}
\psi _{1/2}(t) &=&\frac{1}{\sqrt{\pi t}}\int_{0}^{+\infty }e^{-\lambda y}e^{-%
\frac{y^{2}}{4t}}dy=\Pr \left\{ |B(t)|<U\right\}  \label{pr3} \\
&=&\Pr \left\{ \sup_{0\leq s\leq t}B(s)<U\right\} .  \notag
\end{eqnarray}%
The previous expression can be checked directly by applying (\ref{qua}):%
\begin{eqnarray}
\psi _{1/2}(t) &=&E_{1/2,1}(-\lambda \sqrt{t})  \label{prr} \\
&=&[\text{by the duplication property of the Gamma]}  \notag \\
&=&\sum_{j=0}^{\infty }\frac{\left( -2\lambda \sqrt{t}\right) ^{j}\Gamma
\left( \frac{j}{2}+\frac{1}{2}\right) }{\Gamma (j+1)\sqrt{\pi }}  \notag \\
&=&\frac{1}{\sqrt{\pi }}\int_{0}^{\infty }e^{-z}z^{-\frac{1}{2}%
}\sum_{j=0}^{\infty }\frac{\left( -2\lambda \sqrt{zt}\right) ^{j}}{j!}
\notag \\
&=&\frac{1}{\sqrt{\pi }}\int_{0}^{\infty }e^{-z}z^{-\frac{1}{2}}e^{-2\lambda
\sqrt{zt}},  \notag
\end{eqnarray}%
which gives (\ref{pr3}), after a change of variable.

\

Also for $\nu =1/2^{n},$ $n\geq 1$, the solution can be expressed in terms
of boundary crossing probability of known processes. Indeed the random
process $\mathcal{T}_{\nu }$ coincides in this case with the $(n-1)$-times
\textit{iterated reflecting Brownian motion }defined as $%
I_{n-1}(t)=|B_{1}(|B_{2}(...(|B_{n}(t)|)...)|)|$, where $B_{j}(t)$ are
independent Brownian motions with variance $2t$, for any $j$. The transition
density $q_{1/2^{n}}(y,t)$ of $I_{n-1}$ is given by%
\begin{equation*}
q_{1/2^{n}}(y,t)=\int_{0}^{+\infty }...\int_{0}^{+\infty }\frac{e^{-\frac{%
y^{2}}{4s_{1}}}}{\sqrt{\pi s_{1}}}\frac{e^{-\frac{s_{1}^{2}}{4s_{2}}}}{\sqrt{%
\pi s_{2}}}...\frac{e^{-\frac{s_{n-1}^{2}}{4t}}}{\sqrt{\pi t}}%
ds_{1}...ds_{n-1},\quad y,t>0,
\end{equation*}%
which coincides with the folded solution to the following fractional
diffusion equation%
\begin{equation}
\frac{\partial ^{1/2^{n}}q}{\partial t^{1/2^{n}}}=\frac{\partial ^{2}q}{%
\partial y^{2}},\qquad y\in \mathbb{R},t>0,\quad q(y,0)=\delta (y),
\end{equation}%
(see \cite{ors}, for $n=1$ and \cite{ors2}, for $n>1$). Therefore, in this
case, the solution to the fractional relaxation equation can be expressed as
the crossing probability of an exponential boundary by an iterated
reflecting Brownian motion, i.e.%
\begin{equation*}
\psi _{1/2^{n}}(t)=\psi (I_{n-1}(t))=\int_{0}^{+\infty }e^{-\lambda
y}q_{1/2^{n}}(y,t)dy=\Pr \left\{ I_{n-1}(t)<U\right\} .
\end{equation*}

\

For other rational values of the fractional order $\nu $, such as, for
example $\nu =1/3$, the solution can be still represented as boundary
crossing probability, but of less known processes.

For $\nu =1/3$ the random process $\mathcal{T}_{\nu }$ in (\ref{pr}) reduces
to the process $A(t)$, introduced and studied in \cite{ors2}, whose
transition function is given by%
\begin{equation}
q_{1/3}(y,t)=\sqrt[3]{\frac{3^{2}}{t}}Ai\left( \frac{y}{\sqrt[3]{3t}}\right)
,\quad y,t>0  \label{ai5}
\end{equation}%
where
\begin{equation}
Ai(w)=\frac{1}{\pi }\int_{0}^{\infty }\cos \left( aw+\frac{\alpha ^{3}}{3}%
\right) d\alpha ,\quad w\in \mathbb{R}  \label{ai4}
\end{equation}%
is the Airy function. By exploiting the relationship between (\ref{ai4}) and
the modified Bessel function
\begin{equation*}
I_{\nu }(w)=\sum_{k=0}^{\infty }\frac{\left( \frac{w}{2}\right) ^{2k+\nu }}{%
k!\Gamma (k+\nu +1)},\quad w\in \mathbb{R},
\end{equation*}%
i.e.%
\begin{equation*}
Ai(w)=\frac{\sqrt{w}}{3}\left[ I_{-1/3}\left( \frac{2\sqrt{w^{3}}}{3}\right)
-I_{1/3}\left( \frac{2\sqrt{w^{3}}}{3}\right) \right] ,\quad w>0
\end{equation*}%
we can rewrite the transition density (\ref{ai5}) of the process $A(t),t>0$
as
\begin{equation*}
q_{1/3}(y,t)=\sqrt{\frac{y}{3t}}\left[ I_{-1/3}\left( 2\sqrt{\frac{y}{3^{3}t}%
}\right) -I_{1/3}\left( 2\sqrt{\frac{y}{3^{3}t}}\right) \right] ,\quad y,t>0.
\end{equation*}%
Therefore, in this case, the fractional relaxation can be written as%
\begin{equation*}
\psi _{1/3}(t)=\psi (A(t))=\int_{0}^{+\infty }e^{-\lambda
y}q_{1/3}(y,t)dy=\Pr \left\{ A(t)<U\right\} .
\end{equation*}%
It can be worth comparing the asymptotic behavior of the different crossing
probabilities introduced so far. By using the well-known integral
representation of the Mittag-Leffler function%
\begin{equation}
E_{\nu ,\beta }(-ct^{\nu })=\frac{t^{1-\beta }}{\pi }\int_{0}^{+\infty
}r^{\nu -\beta }e^{-rt}\frac{r^{\nu }\sin (\beta \pi )+c\sin ((\beta -\nu
)\pi )}{r^{2\nu }+2r^{\nu }c\cos (\nu \pi )+c^{2}}dr,  \label{pap}
\end{equation}%
we get the following asymptotic behavior of the solution $\psi _{\nu }$:%
\begin{equation}
\psi _{\nu }(t)\simeq \left\{
\begin{array}{l}
1-\frac{\lambda t^{\nu }}{\Gamma (1+\nu )}\qquad 0<t<<1 \\
\frac{1}{\lambda t^{\nu }\Gamma (1-\nu )},\qquad t\rightarrow \infty
\end{array}%
.\right.   \label{asy4}
\end{equation}%
Therefore the boundary crossing probability of Brownian motion exhibits a
power decay, for $t\rightarrow \infty $, of exponent $1/2$, instead of the
usual exponential decay of the standard relaxation $\psi .$ For the $n$-th
times iterated Brownian motion the exponent $1/2^{n}$ of $t$ is smaller than
$1/2$ and decreases as $n$ becomes larger. This is intuitively explained by
the fact that the number of subordinations increases in the definition of
the process $I_{n}$: this strays the fractional relaxation more and more
away from the standard (exponential) behavior, as $n$ increases, and makes
the tail of the relaxation more and more heavy.

For the process $A(t)$ the crossing probability possesses a power decay, for
$t\rightarrow \infty $, with exponent $1/3$ which is between the Brownian
case and the iterated one (for any $n>1$).

\subsection{Exponential boundary crossing probabilities of more general
processes}

We now present some extensions of the previous results, obtained by
considering the crossing probabilities of different kinds of processes. This
corresponds to substituting the random process $\mathcal{T}_{\nu }(t)$ in (%
\ref{pr}) with some other process, linked to the Brownian motion by various
relationships, such as the elastic Brownian motion, the Bessel process (or
its square), the first passage time through a level $t$ by a standard
Brownian motion or its sojourn time on the positive half line.

We start from the latter, which, being a nondecreasing Lévy process, can be
considered as a subordinator. Let $\Gamma _{t}^{+}(t)=meas\left\{
s<t:B(t)>0\right\} $ be the \textit{sojourn time} on the positive half line
of a standard Brownian motion $B,$ then its density $q^{+}(s,t)$ is given by%
\begin{equation}
q^{+}(s,t)=\Pr \left\{ \Gamma _{t}^{+}\in ds\right\} =\frac{ds}{\pi \sqrt{%
s(t-s)}},\quad 0<s<t.  \label{p}
\end{equation}

\

\noindent \textbf{Theorem 2.2 }\textit{Let }$U$\textit{\ be a random
boundary exponentially distributed, with parameter }$\lambda >0$\textit{.
Then the crossing probability of }$U$\textit{\ by the random process }$%
\Gamma _{t}^{+}(t)$\textit{\ with transition density }$q^{+}(s,t),$\textit{\
is given by}%
\begin{equation}
\psi ^{+}(t)=\psi (\Gamma ^{+}(t))=\Pr \left\{ \Gamma ^{+}(t)<U\right\}
=e^{-\lambda t/2}I_{0}\left( \frac{\lambda t}{2}\right)   \label{psi}
\end{equation}%
\textit{and (\ref{psi}) solves the following second-order differential
equation}%
\begin{equation}
\frac{d^{2}\psi ^{+}}{dt^{2}}+(\lambda +\frac{1}{t})\frac{d\psi ^{+}}{dt}=-%
\frac{\lambda }{2t}\psi ^{+},\quad \psi ^{+}(0)=1.  \label{psi2}
\end{equation}

\noindent \textbf{Proof \ }We write the crossing probability as%
\begin{eqnarray}
\psi ^{+}(t) &=&\int_{0}^{t}e^{-\lambda s}\frac{ds}{\pi \sqrt{s(t-s)}}
\label{effe} \\
&=&\left[ \text{formula }3.383.1,\text{ p.365 }\cite{GR}\right]  \notag \\
&=&_{1}F_{1}\left( \frac{1}{2};1;-\lambda t\right)  \notag
\end{eqnarray}%
where $_{1}F_{1}\left( \alpha ,\gamma ;x\right) $\ denotes the confluent
hypergeometric function defined as
\begin{equation*}
_{1}F_{1}\left( \alpha ;\gamma ;x\right) =1+\sum_{j=1}^{\infty }\frac{\alpha
(\alpha +1)...(\alpha +j-1)}{\gamma (\gamma +1)...(\gamma +j-1)}\frac{x^{j}}{%
j!},
\end{equation*}%
for $x,\alpha \in \mathbb{C}$ and $\gamma \in \mathbb{C}\backslash \mathbb{Z}%
_{0}^{-}.$

By applying the relationship with the Bessel functions (see formula 9.215.2,
p.1086 \cite{GR}) and, after some computations, we get the final form (\ref%
{psi}). As far as the equation satisfied by (\ref{psi}) is concerned, we
recall that $I_{0}(\lambda x)$ coincides with the solution to the following
equation%
\begin{equation}
\frac{d^{2}}{dx^{2}}I_{0}(\lambda x)+\frac{1}{x}\frac{d}{dx}I_{0}(\lambda
x)=\lambda ^{2}I_{0}(\lambda x),  \label{io}
\end{equation}%
as can be easily checked. Therefore, by the transformation $I_{0}\left(
\frac{\lambda t}{2}\right) =e^{\lambda t/2}\psi ^{+}(t)$, from equation (\ref%
{io}) we get (\ref{psi2}), since%
\begin{eqnarray*}
\frac{d}{dt}I_{0}\left( \frac{\lambda t}{2}\right) &=&\frac{\lambda }{2}%
e^{\lambda t/2}\psi ^{+}(t)+e^{\lambda t/2}\frac{d}{dt}\psi ^{+}(t) \\
\frac{d^{2}}{dt^{2}}I_{0}\left( \frac{\lambda t}{2}\right) &=&\frac{\lambda
^{2}}{4}e^{\lambda t/2}\psi ^{+}(t)+\lambda e^{\lambda t/2}\frac{d}{dt}\psi
^{+}(t)+e^{\lambda t/2}\frac{d^{2}}{dt^{2}}\psi ^{+}(t).
\end{eqnarray*}%
Alternatively we can resort to to the form (\ref{effe}) and exploit the fact
that the confluent hypergeometric function $_{1}F_{1}\left( \alpha ;\gamma
;x\right) $ satisfies the following equation:%
\begin{equation}
x\frac{d^{2}}{dx^{2}}\,_{1}F_{1}+(\gamma -x)\frac{d}{dx}\,_{1}F_{1}=\alpha
\,_{1}F_{1}.
\end{equation}%
By taking into account that%
\begin{eqnarray*}
\frac{d}{dt}\,_{1}F_{1}\left( \frac{1}{2};1;-\lambda t\right) &=&-\lambda
\frac{d}{d(-\lambda t)}\,_{1}F_{1}\left( \frac{1}{2};1;-\lambda t\right) \\
\frac{d^{2}}{dt^{2}}\,_{1}F_{1}\left( \frac{1}{2};1;-\lambda t\right)
&=&\lambda ^{2}\frac{d}{d(-\lambda t)^{2}}\,_{1}F_{1}\left( \frac{1}{2}%
;1;-\lambda t\right) ,
\end{eqnarray*}%
we get again (\ref{psi2}).\hfill $\blacksquare $

\

The asymptotic behavior of $\psi ^{+}(t)$ can be deduced by considering that
$I_{\nu }(x)\simeq (x/2)^{\nu }/\Gamma (\nu +1),$ as $x\rightarrow 0,$ and
that
\begin{equation*}
_{1}F_{1}\left( \alpha ;\gamma ,x\right) \simeq \frac{\Gamma (\gamma )}{%
\Gamma (\alpha )}e^{-i\pi \alpha }x^{-\alpha },\quad Re(x)\rightarrow
-\infty
\end{equation*}%
(see \cite{kil}, p.29), thus obtaining the following expressions%
\begin{equation}
\psi ^{+}(t)\simeq \left\{
\begin{array}{l}
1-\frac{\lambda t}{2}\qquad 0<t<<1 \\
\frac{1}{\sqrt{\lambda \pi t}},\qquad t\rightarrow \infty%
\end{array}%
\right. .  \label{psi3}
\end{equation}%
The limiting behavior of $\psi ^{+}(t)$ is the same of a standard
relaxation, for $t\rightarrow 0$, while coincides with that of $\psi
_{1/2}(t)$, for $t\rightarrow \infty $ (up to multiplicative constants).

\

Another process that can be considered instead of the random time $\mathcal{T%
}_{\nu }(t)$ in (\ref{pr}) is the \textit{first passage time through a level
}$t$ by a standard Brownian motion, denoted as%
\begin{equation*}
T(t)=\inf \left\{ s>0:B(s)=t\right\} .
\end{equation*}%
Therefore, we are interested in the following crossing probability%
\begin{equation}
\psi _{T}(t)=\psi (T(t))=\int_{0}^{\infty }e^{-\lambda s}q_{T}(s,t)ds=\Pr
\left\{ T(t)<U\right\} ,  \label{ti}
\end{equation}%
where the density of $T(t),t>0$ is the well-known stable law of index $1/2$,
i.e.%
\begin{equation*}
q_{T}(s,t)=\frac{te^{-t^{2}/2s}}{\sqrt{2\pi s^{3}}},\qquad s,t>0.
\end{equation*}%
Therefore (\ref{ti}) can be easily evaluated, since the Laplace transform of
the first passage time is well-known:%
\begin{equation}
\psi _{T}(t)=e^{-t\sqrt{2\lambda }}.  \label{ti2}
\end{equation}%
Clearly $\psi _{T}(t)$ satisfies the standard relaxation equation, even if
with a different constant:%
\begin{equation*}
\frac{d\psi _{T}}{dt}=-\sqrt{2\lambda }\psi _{T},\quad \psi _{T}(0)=1.
\end{equation*}%
We remark that time-changing the relaxation $\psi $ by the $1/2$-stable
subordinator $T(t)$ produces again a standard relaxation, while performing
the same operation by the inverse stable subordinator $E(t)$ we get the
fractional relaxation $\psi _{\frac{1}{2}}$ (as mentioned in the
introduction).

If we now consider $n$ independent Brownian motions $B_{j},$ $j=1,...,n$ and
construct by them the $n$-times subordinated process $%
T_{1}(T_{2}(...T_{n}(t)...))$, $t>0$, where $T_{j}=\inf \left\{
s>0:B_{j}(s)=t\right\} ,$ $j=1,...,n$, then its crossing probability can be
evaluated as follows:%
\begin{eqnarray}
\psi _{T}^{n}(t) &=&\Pr \left\{ T_{1}(T_{2}(...T_{n}(t)...))<U\right\}
\label{ti3} \\
&=&\int_{0}^{\infty }e^{-\lambda s}\left( \int_{0}^{+\infty
}dz_{1}...\int_{0}^{+\infty }dz_{n-1}\frac{te^{-t^{2}/2z_{1}}}{\sqrt{2\pi
z_{1}^{3}}}...\frac{z_{n-1}e^{-z_{n-1}^{2}/2z_{n}}}{\sqrt{2\pi z_{n}^{3}}}%
\frac{z_{n}e^{-z_{n}^{2}/2s}}{\sqrt{2\pi s^{3}}}\right) ds  \notag \\
&=&\int_{0}^{+\infty }dz_{1}...\int_{0}^{+\infty }dz_{n-1}\frac{%
te^{-t^{2}/2z_{1}}}{\sqrt{2\pi z_{1}^{3}}}...\frac{%
z_{n-1}e^{-z_{n-1}^{2}/2z_{n}}}{\sqrt{2\pi z_{n}^{3}}}\int_{0}^{\infty
}e^{-\lambda s}\frac{z_{n}e^{-z_{n}^{2}/2s}}{\sqrt{2\pi s^{3}}}ds  \notag \\
&=&\int_{0}^{+\infty }dz_{1}...\frac{te^{-t^{2}/2z_{1}}}{\sqrt{2\pi z_{1}^{3}%
}}...\int_{0}^{+\infty }\frac{z_{n-1}e^{-z_{n-1}^{2}/2z_{n}}}{\sqrt{2\pi
z_{n}^{3}}}e^{-z_{n}\sqrt{2\lambda }}dz_{n-1}  \notag \\
&=&\int_{0}^{+\infty }dz_{1}...\frac{te^{-t^{2}/2z_{1}}}{\sqrt{2\pi z_{1}^{3}%
}}...\int_{0}^{+\infty }\frac{z_{n-2}e^{-z_{n-2}^{2}/2z_{n-1}}}{\sqrt{2\pi
z_{n-1}^{3}}}e^{-z_{n-1}\sqrt{2\sqrt{2\lambda }}}dz_{n-2}  \notag \\
&=&e^{-\lambda ^{\frac{1}{2^{n}}}2^{1-\frac{1}{2^{n}}}t}.  \notag
\end{eqnarray}%
Again the probability $\psi _{T}^{n}$ satisfies (for any $n$) the standard
relaxation equation with the constant $\lambda ^{\frac{1}{2^{n}}}2^{1-\frac{1%
}{2^{n}}}$ and displays an asymptotic behavior similar to the standard
relaxation, despite the complicated construction via the $n$-times
subordination.

\

We analyze now the crossing probability of an exponential boundary $U$ by a
\textit{squared Bessel process}. Let us denote by $R_{\gamma }^{2}(t)=\left(
R_{\gamma }(t)\right) ^{2},$ $t>0$ the square of a $\gamma $-Bessel process,
starting at zero. It is well known that, for $\gamma =n$, this process can
be expressed as
\begin{equation*}
R_{n}^{2}(t)=\sum_{j=1}^{n}B_{j}^{2}(t),\qquad t>0,
\end{equation*}%
where $B_{j}(t),j=1,...n,$ are independent Brownian motion in $\mathbb{R}%
^{n}.$ Moreover the density of $R_{\gamma }^{2}\ $can be written as%
\begin{equation*}
p_{\gamma }^{2}(s,t)=\frac{s^{\frac{\gamma }{2}-1}e^{-\frac{s}{2t}}}{(2t)^{%
\frac{\gamma }{2}}\Gamma \left( \frac{\gamma }{2}\right) },\qquad s,t>0
\end{equation*}%
(see, for example, \cite{dov}), which is a more tractable form (for our
aims) than that of $R_{\gamma }.$ Thus the crossing probability of this
process can be easily evaluated as follows:%
\begin{eqnarray}
\psi _{\gamma }(t) &=&\Pr \left\{ R_{\gamma }^{2}(t)<U\right\}  \label{tre.1}
\\
&=&\int_{0}^{\infty }e^{-\lambda s}\frac{s^{\frac{\gamma }{2}-1}e^{-\frac{s}{%
2t}}}{(2t)^{\frac{\gamma }{2}}\Gamma \left( \frac{\gamma }{2}\right) }ds
\notag \\
&=&\frac{1}{\left( 2\lambda t+1\right) ^{\frac{\gamma }{2}}},  \notag
\end{eqnarray}%
which satisfies the following first-order differential equation%
\begin{equation*}
\frac{d}{dt}\psi _{\gamma }=\frac{\gamma \lambda }{2\lambda t+1}\psi
_{\gamma },\qquad \psi _{\gamma }(0)=1.
\end{equation*}

In this case, the behavior of $\psi _{\gamma }(t),$ for increasing (but
still finite) values of $t,$ can be represented as $\psi _{\gamma }(t)\simeq
(k/t)^{\gamma /2}$ (for some constant $k$ and for $0<\gamma <2)$ and thus it
coincides with the one described as \textquotedblleft algebraic
decay\textquotedblright\ and displayed by relaxation processes in complex
material (see, for example, \cite{sci}). On the contrary, for the other
fractional relaxations, this is true only in the limit, for $t\rightarrow
\infty .$ Indeed the function (\ref{tre.1}) coincides with the so-called
Nutting law, which is commonly used to fit the experimental data for the
materials featuring non-standard (i.e. non-Debye) relaxation (see \cite{met}
and the references therein).

\

As we have seen, the generalizations analyzed so far in this section are not
linked to fractional equations; on the other hand, in the following case, we
consider crossing probabilities governed again by fractional equations. Let $%
B^{\alpha }(t),t>0$ be the so called \textit{elastic Brownian motion} with
absorbing rate $\alpha >0$ (see \cite{ITO} and \cite{begors}), defined as%
\begin{equation}
B_{\alpha }^{el}(t)=\left\{
\begin{array}{c}
|B(t)|,\qquad t<T_{\alpha } \\
0,\qquad t\geq T_{\alpha }%
\end{array}%
\right. ,  \label{uno6}
\end{equation}%
where $T_{\alpha }$ is a random time with distribution
\begin{equation}
\Pr \left\{ T_{\alpha }>t|\mathcal{B}_{t}\right\} =e^{-\alpha L(0,t)},\qquad
\alpha >0,  \label{uno7}
\end{equation}%
$\mathcal{B}_{t}=\sigma \left\{ B(s),s\leq t\right\} $ is the natural
filtration and $L(0,t)=\lim {}_{\varepsilon \downarrow 0}\frac{1}{%
2\varepsilon }meas\left\{ s\leq t:|B(t)|<\varepsilon \right\} $ is the local
time in the origin of $B.$ It is well known that its distribution can be
expressed as%
\begin{equation}
q_{\alpha }^{el}(s,t)=2e^{\alpha s}\int_{s}^{+\infty }we^{-\alpha w}\frac{%
e^{-\frac{w^{2}}{2t}}}{\sqrt{2\pi t^{3}}}dw+q_{\alpha }(t)\delta (s),\quad
s,t>0  \label{new2}
\end{equation}%
where $\delta (s)$ is the Dirac's Delta function with pole in the origin and
\begin{equation*}
q_{\alpha }(t)=1-\Pr \left\{ B_{\alpha }^{el}(t)>0\right\} =1-2e^{\frac{%
\alpha ^{2}t}{2}}\int_{\alpha \sqrt{t}}^{+\infty }\frac{e^{-\frac{w^{2}}{2}}%
}{\sqrt{2\pi }}dw
\end{equation*}%
is the probability that the process is absorbed by the barrier in zero up to
time $t.$ Thus we define the crossing probability of an exponential boundary
$U$ by the process $B_{\alpha }^{el}$ as%
\begin{equation}
\psi _{\alpha }^{el}(t)=\Pr \left\{ B_{\alpha }^{el}(t)<U\right\}
=\int_{0}^{\infty }e^{-\lambda s}q_{\alpha }^{el}(s,t)ds.  \label{el}
\end{equation}

\

\noindent \textbf{Theorem 2.3 }\textit{Let }$U$\textit{\ be a random
boundary exponentially distributed, with parameter }$\lambda >0$\textit{.
Then the crossing probability of }$U$\textit{\ by the random process }$%
B_{\alpha }^{el}(t)$\textit{\ with transition density }$q_{\alpha
}^{el}(s,t),$\textit{\ is given, for any }$\lambda \neq \alpha ,$\textit{\ by%
}%
\begin{equation}
\psi _{\alpha }^{el}(t)=\Pr \left\{ B_{\alpha }^{el}(t)<U\right\} =1-\frac{%
\lambda }{\lambda -\alpha }\left[ E_{\frac{1}{2},1}\left( -\frac{\alpha
\sqrt{t}}{\sqrt{2}}\right) -E_{\frac{1}{2},1}\left( -\frac{\lambda \sqrt{t}}{%
\sqrt{2}}\right) \right] ,  \label{el8}
\end{equation}%
\textit{while, for }$\alpha =\lambda ,$\textit{\ it coincides with}%
\begin{equation}
\psi _{\lambda }^{el}(t)=\Pr \left\{ B_{\lambda }^{el}(t)<U\right\}
=1-\lambda \sqrt{2t}E_{\frac{1}{2},\frac{1}{2}}\left( -\frac{\lambda \sqrt{t}%
}{\sqrt{2}}\right) .  \label{el3}
\end{equation}%
\textit{The crossing probability }$\psi _{\alpha }^{el}(t)$\textit{\
satisfies, for any }$\alpha ,\lambda >0,$\textit{\ the following fractional
differential equation}%
\begin{equation}
\frac{d}{dt}\psi _{\alpha }^{el}+\frac{\alpha +\lambda }{\sqrt{2}}\frac{%
d^{1/2}}{dt^{1/2}}\psi _{\alpha }^{el}=\frac{\alpha \lambda }{2}(1-\psi
_{\alpha }^{el})-\frac{\lambda }{\sqrt{2\pi t}},\quad \psi _{\alpha
}^{el}(0)=1.  \label{el9}
\end{equation}%
\textbf{Proof }\ We take the Laplace transform of (\ref{el}), which reads,
for any $\alpha ,\lambda >0$:%
\begin{eqnarray}
&&\int_{0}^{\infty }e^{-\eta t}\psi _{\alpha }^{el}(t)dt=\int_{0}^{\infty
}e^{-\eta t}dt\int_{0}^{\infty }e^{-\lambda s}q_{\alpha }^{el}(s,t)ds
\label{el10} \\
&=&2\int_{0}^{\infty }e^{-\eta t}dt\int_{0}^{\infty }e^{-\lambda s+\alpha
s}ds\int_{s}^{+\infty }we^{-\alpha w}\frac{e^{-\frac{w^{2}}{2t}}}{\sqrt{2\pi
t^{3}}}dw+  \notag \\
&&+\frac{1}{\eta }-2\int_{0}^{\infty }e^{-\eta t+\frac{\alpha ^{2}t}{2}%
}dt\int_{\alpha \sqrt{t}}^{+\infty }\frac{e^{-\frac{w^{2}}{2}}}{\sqrt{2\pi }}%
dw  \notag \\
&=&2\int_{0}^{\infty }e^{-\lambda s+\alpha s}ds\int_{s}^{+\infty
}e^{-(\alpha +\sqrt{2\eta })w}dw+\frac{1}{\eta }-  \notag \\
&&-\frac{2}{2\eta -\alpha ^{2}}+\frac{2\alpha }{\sqrt{2\pi }(2\eta -\alpha
^{2})}\frac{1}{\sqrt{\eta }}\int_{0}^{+\infty }e^{-z}\frac{1}{\sqrt{z}}dz
\notag \\
&=&\frac{2}{\sqrt{2\eta }+\alpha }\int_{0}^{\infty }e^{-\lambda s-\sqrt{%
2\eta }s}ds+\frac{2\eta -\alpha ^{2}-2\eta +\sqrt{2\eta }\alpha }{\eta
(2\eta -\alpha ^{2})}  \notag \\
&=&\frac{2}{(\sqrt{2\eta }+\alpha )(\sqrt{2\eta }+\lambda )}+\frac{\alpha (%
\sqrt{2\eta }-\alpha )}{\eta (2\eta -\alpha ^{2})}  \notag \\
&=&\frac{\alpha \lambda \eta ^{-1}+\sqrt{2}\alpha \eta ^{-\frac{1}{2}}+2}{(%
\sqrt{2\eta }+\alpha )(\sqrt{2\eta }+\lambda )}.  \notag
\end{eqnarray}%
We can check that (\ref{el10}) coincides with the Laplace transform of (\ref%
{el8}), for $\alpha \neq \lambda ,$ as follows:%
\begin{eqnarray*}
&&\mathcal{L}\left\{ \psi _{\alpha }^{el};\eta \right\} =\int_{0}^{\infty
}e^{-\eta t}\psi _{\alpha }^{el}(t)dt \\
&=&\frac{1}{\eta }-\frac{\lambda }{\lambda -\alpha }\sum_{j=0}^{\infty }%
\frac{1}{\Gamma \left( \frac{j}{2}+1\right) }\left[ \left( -\frac{\alpha }{%
\sqrt{2}}\right) ^{j}-\left( -\frac{\lambda }{\sqrt{2}}\right) ^{j}\right]
\int_{0}^{\infty }e^{-\eta t}t^{\frac{j}{2}}dt \\
&=&\frac{1}{\eta }-\frac{\lambda }{\lambda -\alpha }\frac{1}{\eta }%
\sum_{j=0}^{\infty }\left[ \left( -\frac{\alpha }{\sqrt{2\eta }}\right)
^{j}-\left( -\frac{\lambda }{\sqrt{2\eta }}\right) ^{j}\right] \\
&=&\frac{1}{\eta }-\frac{\lambda }{\lambda -\alpha }\frac{1}{\eta }\left[
\frac{\sqrt{2\eta }}{\sqrt{2\eta }+\alpha }-\frac{\sqrt{2\eta }}{\sqrt{2\eta
}+\lambda }\right] ,
\end{eqnarray*}%
which easily gives (\ref{el10}). As a further check of (\ref{el8}), it is
easy to see that, for $\alpha =0$ (in the case of no absorption) it reduces
to $\psi _{\frac{1}{2}}(t)=E_{1/2,1}(-\lambda \sqrt{t})$, since in this case
$B^{el}(t)=|B(t)|,$ $t>0.$

For $\alpha =\lambda $ the Laplace transform (\ref{el10}) becomes%
\begin{equation}
\int_{0}^{\infty }e^{-\eta t}\psi _{\lambda }^{el}(t)dt=\frac{\lambda
^{2}\eta ^{-1}+\sqrt{2}\lambda \eta ^{-\frac{1}{2}}+2}{(\sqrt{2\eta }%
+\lambda )^{2}}.  \label{ml}
\end{equation}%
By comparing (\ref{ml}) with the formula holding for the Laplace transform
of the GML function defined in (\ref{gml2}) (see \cite{kil}, p.47), i.e.%
\begin{equation}
\mathcal{L}\left\{ t^{\gamma -1}E_{\beta ,\gamma }^{\delta }(\omega t^{\beta
});\eta \right\} =\frac{\eta ^{\beta \delta -\gamma }}{(\eta ^{\beta
}-\omega )^{\delta }},  \label{pra}
\end{equation}%
(where $Re(\beta )>0,$ $Re(\gamma )>0,$ $Re(\delta )>0$
and $\eta >|\omega |^{\frac{1}{Re(\beta )}}),$ we easily obtain
\begin{equation}
\psi _{\lambda }^{el}(t)=1-\frac{\lambda \sqrt{t}}{\sqrt{2}}E_{\frac{1}{2},%
\frac{3}{2}}^{2}\left( -\frac{\lambda \sqrt{t}}{\sqrt{2}}\right) ,
\label{ppp}
\end{equation}%
which can be also rewritten as (\ref{el3}).

By taking the Laplace transform of equation (\ref{el8}) and considering the
well-known expression for the Laplace transform of the Caputo derivative,
i.e.%
\begin{eqnarray}
\mathcal{L}\left\{ \frac{d^{\nu }u}{dt^{\nu }};\eta \right\}
&=&\int_{0}^{\infty }e^{-\eta t}\frac{d^{\nu }}{dt^{\nu }}u(t)dt  \label{jen}
\\
&=&\eta ^{\nu }\mathcal{L}\left\{ u;\eta \right\} -\sum_{r=0}^{m-1}\eta
^{\nu -r-1}\left. \frac{d^{r}}{dt^{r}}u(t)\right\vert _{t=0},  \notag
\end{eqnarray}%
we get
\begin{eqnarray}
&&\eta \mathcal{L}\left\{ \psi _{\alpha }^{el};\eta \right\} -\psi _{\alpha
}^{el}(0)+\frac{\alpha +\lambda }{\sqrt{2}}\eta ^{\frac{1}{2}}\mathcal{L}%
\left\{ \psi _{\alpha }^{el};\eta \right\} -\frac{\alpha +\lambda }{\sqrt{2}}%
\eta ^{-\frac{1}{2}}\psi _{\alpha }^{el}(0)  \label{jen2} \\
&=&\frac{\alpha \lambda }{2}(\frac{1}{\eta }-\mathcal{L}\left\{ \psi
_{\alpha }^{el};\eta \right\} )-\frac{\lambda \Gamma \left( \frac{1}{2}%
\right) }{\sqrt{2\pi \eta }}.  \notag
\end{eqnarray}%
By taking account the initial condition $\psi _{\alpha }^{el}(0)=1,$ the
solution of (\ref{jen2}) coincides with (\ref{el10}).\hfill $\blacksquare $

\

In order to study the asymptotics of the solution $\psi _{\alpha }^{el}(t),$
for $\alpha \neq \lambda $, we use the integral expansion for the
Mittag-Leffler function (\ref{pap}), so that we get%
\begin{equation}
\psi _{\alpha }^{el}(t)=1-\frac{\lambda }{\lambda -\alpha }\frac{1}{\pi }%
\int_{0}^{+\infty }z^{-1/2}e^{-z}\left[ \frac{\frac{\alpha }{\sqrt{2}}}{%
\frac{z}{\sqrt{t}}+\frac{\alpha ^{2}}{2}\sqrt{t}}-\frac{\frac{\lambda }{%
\sqrt{2}}}{\frac{z}{\sqrt{t}}+\frac{\lambda ^{2}}{2}\sqrt{t}}\right] dz.
\label{pa}
\end{equation}%
Therefore the limiting behavior of the crossing probability reads%
\begin{equation}
\psi _{\alpha }^{el}(t)\simeq \left\{
\begin{array}{l}
1-\frac{\lambda \sqrt{2t}}{\sqrt{\pi }},\;\qquad 0<t<<1 \\
1-\frac{\sqrt{2}}{\alpha \sqrt{\pi t}},\;\qquad t\rightarrow \infty%
\end{array}%
\right.  \label{pap1}
\end{equation}%
where the first line is obtained from (\ref{pa}) by the following
calculations:%
\begin{eqnarray*}
\psi _{\alpha }^{el}(t) &=&1+\frac{\lambda \sqrt{t}}{\sqrt{2}\pi }%
\int_{0}^{+\infty }z^{-3/2}e^{-z}dz \\
&=&1+\frac{\lambda \sqrt{t}}{\sqrt{2}\pi }\Gamma \left( -\frac{1}{2}\right)
\\
&=&\left[ \text{by the reflection formula of Gamma function}\right] \\
&=&1-\frac{\lambda \sqrt{2t}}{\sqrt{\pi }}.
\end{eqnarray*}%
Thus, in this case, the crossing probability maintains a limiting behavior
similar to the previous ones for $t\rightarrow 0$, but drastically different
for $t\rightarrow \infty $ (see (\ref{asy4})). In the last case instead of
tending to zero, it tends to one: this can be intuitively explained by
considering that the absorbing effect is stronger as $t$ increases and, in
the limit, the process $B^{el}$ will be absorbed with probability one. This
effect is directly correlated with the absorbing rate $\alpha .$ Thus it is
evident from (\ref{pap1}) that $\psi _{\alpha }^{el}$ looses the usual
property of complete monotonicity that characterizes the standard and also
the fractional relaxations (see, for example, \cite{mai}).

In the case $\alpha =\lambda $ we must apply the integral expansion of GML
functions (see \cite{beg3})%
\begin{equation}
E_{\nu ,\beta }^{k}(-ct^{\nu })=\frac{t^{1-\beta }}{2\pi i}\int_{0}^{\infty
}e^{-rt}r^{\nu k-\beta }\left[ \frac{e^{i\pi \beta }}{(r^{\nu }+ce^{i\pi \nu
})^{k}}-\frac{e^{-i\pi \beta }}{(r^{\nu }+ce^{-i\pi \nu })^{k}}\right] dr,
\label{asy}
\end{equation}%
(for $k=2$, $\nu =1/2$, $\beta =3/2$ and $c=\lambda /\sqrt{2})$ so that
formula (\ref{ppp}) can be developed as
\begin{eqnarray*}
&&\psi _{\lambda }^{el}(t) \\
&=&1+\frac{\lambda }{\sqrt{2}}\frac{1}{2\pi }\int_{0}^{\infty }\frac{%
e^{-rt}r^{-\frac{1}{2}}}{(r+\frac{\lambda ^{2}}{2})^{2}}\left[ \left( \sqrt{r%
}-\frac{i\lambda }{\sqrt{2}}\right) ^{2}+\left( \sqrt{r}+\frac{i\lambda }{%
\sqrt{2}}\right) ^{2}\right] dr \\
&=&1+\frac{\lambda }{\sqrt{2t}}\frac{1}{\pi }\int_{0}^{\infty }e^{-z}z^{-%
\frac{1}{2}}\frac{\frac{z}{t}-\frac{\lambda ^{2}}{2}}{(\frac{z}{t}+\frac{%
\lambda ^{2}}{2})^{2}}dz.
\end{eqnarray*}%
Therefore, also for $\alpha =\lambda $, the asymptotic behavior is given
exactly by (\ref{pap1}).

\

\noindent \textbf{Remark 2.1 }An interesting link can be found between the
crossing probabilities $\psi _{\alpha }^{el}(t)$ and $\psi _{1/2}(t)$: for $%
\lambda \neq \alpha ,$ the first one can be rewritten, in view of (\ref{el8}%
) and (\ref{prr}), as%
\begin{equation}
\psi _{\alpha }^{el}(t)=1-\frac{\lambda }{\lambda -\alpha }\left[ \psi
_{1/2}^{\alpha }(t)-\psi _{1/2}^{\lambda }(t)\right] ,
\end{equation}%
where $\psi _{1/2}^{\alpha }(t)$ and $\psi _{1/2}^{\lambda }(t)$ indicate
the crossing probability $\Pr \left\{ |B(t)|<U\right\} $ of an exponential
boundary $U$ of parameter $\alpha $ and $\lambda $, respectively, by a
Brownian motion$.$ Thus the following identity is also verified for the
corresponding differential equations:%
\begin{eqnarray*}
\frac{d^{1/2}}{dt^{1/2}}\psi _{\alpha }^{el} &=&-\frac{\lambda }{\lambda
-\alpha }\left[ \frac{d^{1/2}}{dt^{1/2}}\psi _{1/2}^{\alpha }(t)-\frac{%
d^{1/2}}{dt^{1/2}}\psi _{1/2}^{\lambda }(t)\right] \\
&=&\frac{\lambda }{\lambda -\alpha }\left[ \frac{\alpha }{\sqrt{2}}\psi
_{1/2}^{\alpha }(t)-\frac{\lambda }{\sqrt{2}}\psi _{1/2}^{\lambda }(t)\right]
,
\end{eqnarray*}%
by applying Theorem 2.1, for $\nu =1/2$.

\subsection{Crossing probabilities of a Gamma distributed boundary}

We extend the previous results by considering the crossing probabilities of
a random boundary, distributed with different laws, instead of the
exponential one. In particular we choose its natural generalization, i.e.
the Gamma distribution. Thus we are considering the following probability,
which extends formula (\ref{ott})%
\begin{equation}
\psi _{\frac{1}{2}}^{k}(t)=\Pr \left\{ |B(t)|<G\right\} =\int_{0}^{\infty }%
\left[ 1-F_{G}(y)\right] \frac{e^{-y^{2}/4t}}{\sqrt{\pi t}}dy,  \label{g1}
\end{equation}%
where $G$ is a Gamma r.v. with parameters $\lambda ,k>0$ and $F_{G}$ denotes
its cumulative distribution function. For our convenience, we write the
latter as follows:%
\begin{equation}
F_{G}(y)=\frac{\lambda ^{k}}{\Gamma (k)}\int_{0}^{y}e^{-\lambda z}z^{k-1}dz=%
\frac{(\lambda y)^{k}}{\Gamma (k)}\sum_{j=0}^{\infty }\frac{(-\lambda y)^{j}%
}{j!(j+k)}.  \label{g2}
\end{equation}

\

\noindent \textbf{Theorem 2.4}\ \textit{Let }$G$\textit{\ be a random
boundary distributed as a Gamma with parameters }$\lambda ,k>0.$\textit{\
Then the crossing probability of }$G$\textit{\ by a standard Brownian motion
is given by}%
\begin{equation}
\psi _{\frac{1}{2}}^{k}(t)=\Pr \left\{ |B(t)|<G\right\} =1-(\lambda \sqrt{t}%
)^{k}E_{\frac{1}{2},\frac{k}{2}+1}^{k}(-\lambda \sqrt{t}),  \label{g5}
\end{equation}%
\textit{which satisfies the following fractional relaxation equation}%
\begin{equation}
\sum_{j=1}^{k}\binom{k}{j}\lambda ^{-j}\frac{d^{\frac{j}{2}}}{dt^{\frac{j}{2}%
}}\psi _{\frac{1}{2}}^{k}(t)=-\psi _{\frac{1}{2}}^{k}(t),  \label{ma}
\end{equation}%
\textit{with initial condition }$\psi _{\frac{1}{2}}^{k}(0)=1,$\textit{\ for
}$k\geq 1$\textit{, and the additional conditions }%
\begin{eqnarray*}
\left. \frac{d^{r}}{dt^{r}}\psi _{\frac{1}{2}}^{k}(t)\right\vert _{t=0} &=&0,%
\mathit{\ }r=1,...,\left\lfloor \frac{k}{2}\right\rfloor ,\mathit{\ }\text{%
\textit{for any odd}}\mathit{\ }k>1 \\
\left. \frac{d^{r}}{dt^{r}}\psi _{\frac{1}{2}}^{k}(t)\right\vert _{t=0} &=&0,%
\mathit{\ }r=1,...,\frac{k}{2}-1,\mathit{\ }\text{\textit{for any even}}%
\mathit{\ }k>2.
\end{eqnarray*}

\noindent \textbf{Proof} We can rewrite (\ref{g1}) as%
\begin{eqnarray}
\psi _{\frac{1}{2}}^{k}(t) &=&\int_{0}^{\infty }\left[ 1-\frac{(\lambda
y)^{k}}{\Gamma (k)}\sum_{j=0}^{\infty }\frac{(-\lambda y)^{j}}{j!(j+k)}%
\right] \frac{e^{-y^{2}/4t}}{\sqrt{\pi t}}dy  \label{g3} \\
&=&1-\frac{1}{\Gamma (k)\sqrt{\pi t}}\sum_{j=0}^{\infty }\frac{%
(-1)^{j}\lambda ^{j+k}}{j!(j+k)}\int_{0}^{\infty }y^{j+k}e^{-y^{2}/4t}dy
\notag \\
&=&1-\frac{1}{\Gamma (k)\sqrt{\pi }}\sum_{j=0}^{\infty }\frac{%
(-1)^{j}(2\lambda \sqrt{t})^{j+k}}{j!(j+k)}\Gamma \left( \frac{j}{2}+\frac{k%
}{2}+\frac{1}{2}\right)  \notag \\
&=&1-\frac{2}{\Gamma (k)}\sum_{j=0}^{\infty }\frac{(-1)^{j}(\lambda \sqrt{t}%
)^{j+k}}{j!(j+k)}\frac{\Gamma \left( j+k\right) }{\Gamma \left( \frac{j}{2}+%
\frac{k}{2}\right) }  \notag \\
&=&1-\frac{(\lambda \sqrt{t})^{k}}{\Gamma (k)}\sum_{j=0}^{\infty }\frac{%
\Gamma \left( j+k\right) (-\lambda \sqrt{t})^{j}}{j!\Gamma \left( \frac{j}{2}%
+\frac{k}{2}+1\right) }.  \notag
\end{eqnarray}%
If we now assume that $k$ is an integer, we can recognize in (\ref{g3}) the
GML function (\ref{gml2}), so that we get (\ref{g5}). As a further check, it
is easy to ascertain that, in the special case $k=1$ (where the r.v. $G$
reduces to the exponential r.v. $U$), the crossing probability $\psi _{\frac{%
1}{2}}^{k}$ given in (\ref{g5}) coincides with the fractional relaxation $%
\psi _{\frac{1}{2}}$ in (\ref{prr}):%
\begin{eqnarray}
\psi _{\frac{1}{2}}^{k}(t) &=&1-\lambda \sqrt{t}E_{\frac{1}{2},\frac{3}{2}%
}(-\lambda \sqrt{t}) \\
&=&1+\sum_{l=1}^{\infty }\frac{(-\lambda \sqrt{t})^{l}}{\Gamma \left( \frac{l%
}{2}+1\right) }=E_{\frac{1}{2},1}(-\lambda \sqrt{t})=\psi _{\frac{1}{2}}(t).
\notag
\end{eqnarray}%
In order to derive equation (\ref{ma}) we resort to the Laplace transform of
(\ref{g5}) which reads:%
\begin{equation}
\mathcal{L}\left\{ \psi _{\frac{1}{2}}^{k};\eta \right\} =\frac{(\sqrt{\eta }%
+\lambda )^{k}-\lambda ^{k}}{\eta (\sqrt{\eta }+\lambda )^{k}},  \label{ma2}
\end{equation}%
by applying again formula (\ref{pra}), for $\gamma =\frac{k}{2}+1,$ $\beta =%
\frac{1}{2}$ and $\delta =k.$ We now rewrite (\ref{ma2}) as follows%
\begin{equation}
\sum_{j=0}^{k}\binom{k}{j}\lambda ^{k-j}\left[ \eta ^{\frac{j}{2}}\mathcal{L}%
\left\{ \psi _{\frac{1}{2}}^{k};\eta \right\} -\eta ^{\frac{j}{2}-1}\right]
=-\frac{\lambda ^{k}}{\eta }.  \label{ma5}
\end{equation}%
By simplifying this expression, we can recognize the Laplace transform of
equation (\ref{ma}). We can check that the initial conditions are satisfied,
by using the series expression of $E_{\nu ,\beta }^{k}(-ct^{\nu })$, and
considering that for $t=0,$ $E_{\nu ,\beta }^{k}(-\lambda \sqrt{t})=1/\Gamma
(\beta )$: thus we get
\begin{equation*}
\left. \psi _{\frac{1}{2}}^{k}(t)\right\vert _{t=0}=1-\left. \frac{(\lambda
\sqrt{t})^{k}}{\Gamma \left( \frac{k}{2}+1\right) }\right\vert _{t=0}=1.
\end{equation*}%
For the other conditions, we can apply the following formula of the $r$-th
order derivatives of a GML function (see formula (1.9.6), p.46 of \cite{kil}%
):%
\begin{equation}
\frac{d^{r}}{dz^{r}}\left[ z^{\beta -1}E_{\alpha ,\beta }^{\rho }(\lambda
z^{\alpha })\right] =z^{\beta -r-1}E_{\alpha ,\beta -r}^{\rho }(\lambda
z^{\alpha }),\qquad \lambda \in \mathbb{C},\text{ }r\in \mathbb{N},
\label{derr}
\end{equation}%
so that we get
\begin{equation}
\frac{d^{r}}{dt^{r}}\psi _{\frac{1}{2}}^{k}(t)=-\lambda ^{k}t^{\frac{k}{2}%
-r}E_{\frac{1}{2},\frac{k}{2}-r+1}^{k}(-\lambda \sqrt{t}),\qquad r\in
\mathbb{N}.  \label{der2}
\end{equation}%
By recalling formula (\ref{jen}), we notice that the Laplace form (\ref{ma5}%
) holds if the derivatives of order $r$ of $\psi _{\frac{1}{2}}^{k}$
vanishes for $r=1,...\left\lfloor \frac{k}{2}\right\rfloor $ if $k>1$ is odd
and for $r=1,...\frac{k}{2}-1$ if $k>2$ is even; this is verified by (\ref%
{der2}).

Finally we check that equation (\ref{ma}) becomes, for $k=1,$ the fractional
relaxation equation $\frac{d^{\frac{1}{2}}}{dt^{\frac{1}{2}}}\psi _{\frac{1}{%
2}}(t)=-\lambda \psi _{\frac{1}{2}}(t).$\hfill $\blacksquare $

\ \

\noindent \textbf{Remark 2.2 }By comparing (\ref{g5}) with the results in
\cite{beg2}, we can deduce that the crossing probability $\psi _{\frac{1}{2}%
}^{k}(t)$ can be written in terms of the fractional Poisson process of order
$\nu =\frac{1}{2},$ as%
\begin{equation}
\psi _{\frac{1}{2}}^{k}(t)=\Pr \left\{ T_{k}>t\right\} =\Pr \left\{ \mathcal{%
N}_{\frac{1}{2}}(t)<k\right\} ,  \label{poi}
\end{equation}%
where $T_{k}=\inf \left\{ t>0:\mathcal{N}_{\frac{1}{2}}(t)=k\right\} $ is
the waiting probability of the $k$-th event. On the other hand we can prove
that the following relationship holds between the crossing probabilities
given in (\ref{g1}) for a Gamma boundary of parameters $\left( \lambda
,k\right) $ and $\left( \lambda ,k-1\right) $ (respectively denoted as $\psi
_{\frac{1}{2}}^{k}(t)$ and $\psi _{\frac{1}{2}}^{k-1}(t))$:%
\begin{equation}
\frac{d^{1/2}}{dt^{1/2}}\psi _{\frac{1}{2}}^{k}(t)=-\lambda \left[ \psi _{%
\frac{1}{2}}^{k}(t)-\psi _{\frac{1}{2}}^{k-1}(t)\right] .  \label{poi2}
\end{equation}%
Indeed we can evaluate the fractional derivative of order $1/2$ of $\psi _{%
\frac{1}{2}}^{k}$, by considering (\ref{der2}):%
\begin{eqnarray}
\frac{d^{1/2}}{dt^{1/2}}\psi _{\frac{1}{2}}^{k}(t) &=&-\frac{\lambda ^{k}}{%
\sqrt{\pi }(k-1)!}\sum_{j=0}^{\infty }\frac{(j+k-1)!(-\lambda )^{j}}{%
j!\Gamma \left( \frac{j}{2}+\frac{k}{2}\right) }\int_{0}^{t}(t-s)^{-\frac{1}{%
2}}s^{\frac{k}{2}+\frac{j}{2}-1}ds  \notag \\
&=&-\lambda ^{k}t^{\frac{k}{2}-\frac{1}{2}}E_{\frac{1}{2},\frac{k}{2}+\frac{1%
}{2}}^{k}(-\lambda \sqrt{t}).  \label{g7}
\end{eqnarray}%
By applying to (\ref{g7}) the following recursive formula for GML function
proved in (\cite{beg2})%
\begin{equation}
x^{n}E_{\nu ,n\nu +z}^{m}(-x)+x^{n+1}E_{\nu ,(n+1)\nu
+z}^{m}(-x)=x^{n}E_{\nu ,n\nu +z}^{m-1}(-x),\quad \text{ }n,m>0,z\geq 0,x>0,
\label{gen}
\end{equation}%
for $m=n=k$, $x=-\lambda \sqrt{t},$ $\nu =1/2,$ $z=1/2$, we can rewrite%
\begin{eqnarray}
\frac{d^{1/2}}{dt^{1/2}}\psi _{\frac{1}{2}}^{k}(t) &=&-t^{-\frac{1}{2}%
}\left( \lambda ^{k}t^{\frac{k}{2}}E_{\frac{1}{2},\frac{k}{2}+\frac{1}{2}%
}^{k}(-\lambda \sqrt{t})\right)  \label{g8} \\
&=&-t^{-\frac{1}{2}}\left[ \lambda ^{k}t^{\frac{k}{2}}E_{\frac{1}{2},\frac{k%
}{2}+\frac{1}{2}}^{k-1}(-\lambda \sqrt{t})-\lambda ^{k+1}t^{\frac{k+1}{2}}E_{%
\frac{1}{2},\frac{k}{2}+1}^{k}(-\lambda \sqrt{t})\right]  \notag \\
&=&-\lambda ^{k}t^{\frac{k}{2}-\frac{1}{2}}E_{\frac{1}{2},\frac{k}{2}+\frac{1%
}{2}}^{k-1}(-\lambda \sqrt{t})+\lambda (1-\psi _{\frac{1}{2}}^{k}(t)),
\notag
\end{eqnarray}%
which gives (\ref{poi2}). The latter could be alternatively obtained by
considering that%
\begin{equation*}
p_{k}^{1/2}(t)=\Pr \left\{ \mathcal{N}_{\frac{1}{2}}(t)=k\right\} =\psi _{%
\frac{1}{2}}^{k}(t)-\psi _{\frac{1}{2}}^{k-1}(t)
\end{equation*}%
satisfies (\ref{ai}) with $\nu =1/2$ and taking into account (\ref{poi}).

\

The asymptotic behavior of the crossing probability $\psi _{\frac{1}{2}}^{k}$
for small $t$ can be deduced by the series expression of the GML function%
\begin{equation}
E_{\nu ,\beta }^{k}(-ct^{\nu })\simeq \frac{1}{\Gamma (\beta )}-\frac{%
ct^{\nu }k}{\Gamma (\beta +\nu )},\qquad 0<t<<1,  \label{gen2}
\end{equation}%
so that we get
\begin{equation}
\psi _{\frac{1}{2}}^{k}(t)\simeq 1-\frac{(\lambda \sqrt{t})^{k}}{\Gamma
\left( \frac{k}{2}+1\right) }.  \label{ma4}
\end{equation}%
The same result can be obtained by resorting to the Laplace transform and to
the Tauberian theory, which permits to infer (formally) the asymptotic
behavior of a function $f(t),$ for $t\rightarrow \infty $ and $t\rightarrow
0^{+},$ from the limiting behavior of its Laplace transform $\mathcal{L}%
\left\{ f;\eta \right\} $\ for $\eta \rightarrow 0^{+}$ and $\eta
\rightarrow \infty $, respectively (see also \cite{mai}, for details). To
this aim, we rewrite (\ref{ma2}) as
\begin{equation}
\mathcal{L}\left\{ \psi _{\frac{1}{2}}^{k};\eta \right\} =\frac{1}{\eta }-%
\frac{\lambda ^{k}}{\eta (\sqrt{\eta }+\lambda )^{k}},  \label{ma3}
\end{equation}%
which, for $\eta \rightarrow \infty $, can be approximated as follows%
\begin{equation}
\mathcal{L}\left\{ \psi _{\frac{1}{2}}^{k};\eta \right\} =\frac{1}{\eta }-%
\frac{\lambda ^{k}}{\eta ^{\frac{k}{2}+1}}+o(\eta ^{-\frac{k}{2}-1})
\end{equation}%
so that we get again (\ref{ma4}). For $t\rightarrow \infty $, it is worth
writing (\ref{ma2}) as%
\begin{equation*}
\mathcal{L}\left\{ \psi _{\frac{1}{2}}^{k};\eta \right\} =\frac{%
\sum_{j=1}^{k}\binom{k}{j}\eta ^{\frac{j}{2}-\frac{1}{2}}\lambda ^{-j}}{%
\sum_{j=0}^{k}\binom{k}{j}\eta ^{\frac{j}{2}+\frac{1}{2}}\lambda ^{-j}}%
\simeq \frac{k}{\lambda \eta ^{\frac{1}{2}}},\qquad \eta \rightarrow 0^{+}
\end{equation*}%
so that we get $\psi _{\frac{1}{2}}^{k}(t)\simeq \frac{k}{\lambda \sqrt{\pi t%
}}.$ Thus the limiting behavior of $\psi _{\frac{1}{2}}^{k}$ can be summed
up as follows:%
\begin{equation}
\psi _{\frac{1}{2}}^{k}(t)\simeq \left\{
\begin{array}{l}
1-\frac{(\lambda \sqrt{t})^{k}}{\Gamma \left( \frac{k}{2}+1\right) },\qquad
0<t<<1 \\
\frac{k}{\lambda \sqrt{\pi t}},\qquad t\rightarrow +\infty%
\end{array}%
\right. ,  \label{gen3}
\end{equation}%
which, of course, coincides with (\ref{asy4}) for $k=1$ and $\nu =1/2.$ We
can deduce that, while for small $t$ passing from an exponential boundary to
a Gamma-distributed one makes a relevant difference, for large $t$ this
effect fades away. Indeed the rate of the decreasing to zero for $%
t\rightarrow \infty $ of the crossing probability is exactly the same for
any $k\geq 1.$

\

Analogously, we can generalize the results of Theorem 2.3, by considering
the crossing probability of a Gamma distributed boundary by the elastic
Brownian motion defined in (\ref{uno6}).

\

\noindent \textbf{Theorem 2.5 }\textit{Let }$G$\textit{\ be a random
boundary distributed as a Gamma with parameters }$\lambda ,k>0$\textit{,
then the crossing probability of }$G$\textit{\ by the random process }$%
B_{\alpha }^{el}(t)$\textit{\ with transition density }$q^{el}(s,t)$ \textit{%
(given in (\ref{new2}}))$,$\textit{\ for any }$\lambda ,\alpha >0,$\textit{\
is equal to}%
\begin{equation}
\psi _{k,\alpha }^{el}(t)=\Pr \left\{ B_{\alpha }^{el}(t)<G\right\}
=1-\left( \frac{\lambda \sqrt{t}}{\sqrt{2}}\right) ^{k}\sum_{l=0}^{\infty
}\left( -\frac{\alpha \sqrt{t}}{\sqrt{2}}\right) ^{l}E_{\frac{1}{2},\frac{l+k%
}{2}+1}^{k}\left( -\frac{\lambda \sqrt{t}}{\sqrt{2}}\right) ,  \label{gam1}
\end{equation}%
\textit{which, in the particular case }$\alpha =\lambda ,$\textit{\ reduces
to}%
\begin{equation}
\psi _{k,\lambda }^{el}(t)=\Pr \left\{ B_{\lambda }^{el}(t)<G\right\}
=1-\left( \frac{\lambda \sqrt{t}}{\sqrt{2}}\right) ^{k}E_{\frac{1}{2},\frac{k%
}{2}+1}^{k+1}\left( -\frac{\lambda \sqrt{t}}{\sqrt{2}}\right) .  \label{gam2}
\end{equation}%
\textbf{Proof }\ By following some steps similar to those of Theorem 2.3, we
can write the Laplace transform of $\psi _{k,\alpha }^{el}(t)$ as follows%
\begin{eqnarray}
&&\mathcal{L}\left\{ \psi _{k,\alpha }^{el};\eta \right\} =\int_{0}^{\infty
}e^{-\eta t}dt\int_{0}^{\infty }\left[ 1-F_{G}(s)\right] q_{\alpha
}^{el}(s,t)ds  \label{gam} \\
&=&2\int_{0}^{\infty }\left[ 1-F_{G}(s)\right] e^{\alpha
s}ds\int_{s}^{+\infty }e^{-(\alpha +\sqrt{2\eta })w}dw+\frac{1}{\eta }-
\notag \\
&&-\frac{2}{2\eta -\alpha ^{2}}+\frac{2\alpha }{\sqrt{2\eta }(2\eta -\alpha
^{2})}  \notag \\
&=&\frac{2}{(\sqrt{2\eta }+\alpha )\sqrt{2\eta }}-\frac{2\lambda ^{2}}{\sqrt{%
2\eta }^{k+1}(\sqrt{2\eta }+\alpha )}\sum_{j=0}^{\infty }\binom{k+j-1}{j}%
\left( -\frac{\lambda }{\sqrt{2\eta }}\right) ^{j}+  \notag \\
&&+\frac{\alpha (\sqrt{2\eta }-\alpha )}{\eta (2\eta -\alpha ^{2})}  \notag
\\
&=&\frac{2(\sqrt{2\eta }+\lambda )^{k}-\lambda ^{k}}{\sqrt{2\eta }(\sqrt{%
2\eta }+\alpha )(\sqrt{2\eta }+\lambda )^{k}}+\frac{\alpha }{\eta (\sqrt{%
2\eta }+\alpha )}  \notag \\
&=&\frac{1}{\eta }-\frac{\sqrt{2}\lambda ^{k}}{\sqrt{\eta }(\sqrt{2\eta }%
+\alpha )(\sqrt{2\eta }+\lambda )^{k}}.  \notag
\end{eqnarray}%
We can invert (\ref{gam}) by applying again (\ref{pra}):%
\begin{eqnarray*}
&&\psi _{k,\alpha }^{el}(t) \\
&=&1-\sqrt{2}\lambda ^{k}\mathcal{L}\left\{ \frac{1}{(\sqrt{2\eta }+\alpha )}%
\frac{\eta ^{-1/2}}{(\sqrt{2\eta }+\lambda )^{k}};t\right\} \\
&=&1-\left( \frac{\lambda }{\sqrt{2}}\right) ^{k}\int_{0}^{t}(t-s)^{-1/2}E_{%
\frac{1}{2},\frac{1}{2}}\left( -\frac{\alpha \sqrt{t-s}}{\sqrt{2}}\right) s^{%
\frac{k}{2}-\frac{1}{2}}E_{\frac{1}{2},\frac{k}{2}+\frac{1}{2}}^{k}\left( -%
\frac{\lambda \sqrt{s}}{\sqrt{2}}\right) ds \\
&=&1-\left( \frac{\lambda }{\sqrt{2}}\right) ^{k}\sum_{l=0}^{\infty }\frac{%
\left( -\frac{\alpha }{\sqrt{2}}\right) ^{l}}{\Gamma \left( \frac{l}{2}+%
\frac{1}{2}\right) }\sum_{j=0}^{\infty }\frac{(k+j-1)!\left( -\frac{\lambda
}{\sqrt{2}}\right) ^{j}}{(k-1)!j!\Gamma \left( \frac{j}{2}+\frac{k+1}{2}%
\right) }\times \\
&&\times \int_{0}^{t}(t-s)^{\frac{l}{2}-\frac{1}{2}}s^{\frac{k-1}{2}+\frac{j%
}{2}}ds,
\end{eqnarray*}%
which, after some simplifications, coincides with (\ref{gam1}). For $\alpha
=\lambda $, we can rewrite the latter as follows:%
\begin{eqnarray*}
\psi _{k,\alpha }^{el}(t) &=&1-\left( \frac{\lambda \sqrt{t}}{\sqrt{2}}%
\right) ^{k}\sum_{l=0}^{\infty }\sum_{j=0}^{\infty }\frac{(k+j-1)!\left( -%
\frac{\lambda \sqrt{t}}{\sqrt{2}}\right) ^{j+l}}{(k-1)!j!\Gamma \left( \frac{%
j}{2}+\frac{l+k}{2}+1\right) } \\
&=&1-\left( \frac{\lambda \sqrt{t}}{\sqrt{2}}\right) ^{k}\sum_{l=0}^{\infty
}\sum_{m=l}^{\infty }\frac{(k+m-l-1)!\left( -\frac{\lambda \sqrt{t}}{\sqrt{2}%
}\right) ^{m}}{(k-1)!(m-l)!\Gamma \left( \frac{m}{2}+\frac{k}{2}+1\right) }
\\
&=&1-\left( \frac{\lambda \sqrt{t}}{\sqrt{2}}\right) ^{k}\sum_{m=0}^{\infty }%
\frac{\left( -\frac{\lambda \sqrt{t}}{\sqrt{2}}\right) ^{m}}{\Gamma \left(
\frac{m}{2}+\frac{k}{2}+1\right) }\sum_{l=0}^{m}\binom{k+m-l-1}{m-l} \\
\end{eqnarray*}%
\begin{eqnarray*}
&=&\text{[by the identity proved in \cite{beg4}, p.10]} \\
&=&1-\left( \frac{\lambda \sqrt{t}}{\sqrt{2}}\right) ^{k}\sum_{m=0}^{\infty }%
\frac{\left( -\frac{\lambda \sqrt{t}}{\sqrt{2}}\right) ^{m}}{\Gamma \left(
\frac{m}{2}+\frac{k}{2}+1\right) }\binom{k+m}{k} \\
&=&\psi _{k,\lambda }^{el}(t).
\end{eqnarray*}%
As a final check, we can ascertain that, for $k=1,$ formulae (\ref{gam1})
and (\ref{gam2}) reduce to the corresponding expressions given for the
exponential case in (\ref{el8}) and (\ref{el3}), respectively: indeed (\ref%
{gam1}) can be rewritten, for $k=1,$ as%
\begin{eqnarray*}
\psi _{1,\alpha }^{el}(t) &=&1+\sum_{l=0}^{\infty }\left( -\frac{\alpha
\sqrt{t}}{\sqrt{2}}\right) ^{l}\sum_{j=0}^{\infty }\frac{\left( -\frac{%
\lambda \sqrt{t}}{\sqrt{2}}\right) ^{j+1}}{\Gamma \left( \frac{j+1}{2}+\frac{%
l}{2}+1\right) } \\
&=&1-\sum_{l=0}^{\infty }\frac{\left( -\frac{\alpha \sqrt{t}}{\sqrt{2}}%
\right) ^{l}}{\Gamma \left( \frac{l}{2}+1\right) }+\sum_{l=0}^{\infty
}\left( -\frac{\alpha \sqrt{t}}{\sqrt{2}}\right) ^{l}\sum_{m=0}^{\infty }%
\frac{\left( -\frac{\lambda \sqrt{t}}{\sqrt{2}}\right) ^{m}}{\Gamma \left(
\frac{m+l}{2}+1\right) } \\
&=&1-E_{\frac{1}{2},1}\left( -\frac{\alpha \sqrt{t}}{\sqrt{2}}\right)
+\sum_{l=0}^{\infty }\left( -\frac{\alpha \sqrt{t}}{\sqrt{2}}\right)
^{l}\sum_{k=l}^{\infty }\frac{\left( -\frac{\lambda \sqrt{t}}{\sqrt{2}}%
\right) ^{k-l}}{\Gamma \left( \frac{k}{2}+1\right) } \\
&=&1-E_{\frac{1}{2},1}\left( -\frac{\alpha \sqrt{t}}{\sqrt{2}}\right)
+\sum_{k=0}^{\infty }\frac{\left( -\frac{\lambda \sqrt{t}}{\sqrt{2}}\right)
^{k}}{\Gamma \left( \frac{k}{2}+1\right) }\sum_{l=0}^{k}\left( \frac{\alpha
}{\lambda }\right) ^{l},
\end{eqnarray*}%
which coincides with (\ref{el8}). Formula (\ref{gam2}) immediately reduces
to the expression (\ref{ppp}), for $k=1.$

Finally, putting $\alpha =0$ and substituting $\lambda /\sqrt{2}$ with $%
\lambda ,$ formula (\ref{gam1}) coincides with the corresponding crossing
probability (\ref{g5}), which has been obtained in the case of a free
Brownian motion (with no absorption).$\hfill \blacksquare $

\

The asymptotic behavior of $\psi _{k,\lambda }^{el},$ for small $t$, can be
derived from (\ref{gam2}), by applying again formula (\ref{gen2}).
Alternatively we can use the Laplace transform (\ref{gam}), which can be
approximated as follows, for $\eta \rightarrow \infty $%
\begin{equation}
\mathcal{L}\left\{ \psi _{k,\alpha }^{el};\eta \right\} \simeq \frac{1}{\eta
}-\frac{\lambda ^{k}}{2^{\frac{k}{2}}\eta ^{\frac{k}{2}+1}}.  \notag
\end{equation}%
In both ways, we get the first line of the following formula:%
\begin{equation}
\psi _{k,\lambda }^{el}(t)\simeq \left\{
\begin{array}{l}
1-\left( \frac{\lambda \sqrt{t}}{\sqrt{2}}\right) ^{k}\frac{1}{\Gamma \left(
\frac{k}{2}+1\right) },\qquad 0<t<<1 \\
1-\frac{\sqrt{2}}{\alpha \sqrt{\pi t}},\qquad t\rightarrow +\infty%
\end{array}%
\right. ,  \label{gam6}
\end{equation}%
The second line of the previous expression has been obtained from (\ref{gam}%
), which can be rewritten as
\begin{eqnarray*}
\mathcal{L}\left\{ \psi _{k,\alpha }^{el};\eta \right\} &=&\frac{1}{\eta }-%
\frac{\sqrt{2}}{\sqrt{\eta }\left[ \sum_{j=0}^{k}\binom{k}{j}(2\eta )^{\frac{%
j}{2}+\frac{1}{2}}\lambda ^{-j}+\alpha \sum_{j=0}^{k}\binom{k}{j}(2\eta )^{%
\frac{j}{2}}\lambda ^{-j}\right] } \\
&\simeq &\frac{1}{\eta }-\frac{\sqrt{2}}{\alpha \sqrt{\eta }},\qquad \eta
\rightarrow 0^{+}.
\end{eqnarray*}%
For $k=1,$ formula (\ref{gam6}) coincides with (\ref{pap1}), as was
expected. We finally note that, also in this case, as for the Brownian
motion, the leading term in the expression obtained for $t\rightarrow \infty
$ does not depend on $k$ and thus, for large values of $t$, considering an
exponential or a Gamma distributed boundary does not entail any consequence.

\

The fractional equations satisfied by the crossing probabilities obtained
above can be derived by properly rewriting the Laplace transform in (\ref%
{gam}), as the following theorem shows.

\

\noindent \textbf{Theorem 2.6 }\textit{The crossing probability }$\psi
_{k,\alpha }^{el}$\textit{\ given in (\ref{gam1}) satisfies, for any }$%
\lambda ,\alpha >0,$\textit{\ the following fractional equation}%
\begin{equation}
\sum_{j=0}^{k}\binom{k}{j}\left( \frac{\sqrt{2}}{\lambda }\right) ^{j}\frac{%
d^{\frac{j}{2}+\frac{1}{2}}}{dt^{\frac{j}{2}+\frac{1}{2}}}\psi _{k,\alpha
}^{el}+\frac{\alpha }{\sqrt{2}}\sum_{j=1}^{k}\binom{k}{j}\left( \frac{\sqrt{2%
}}{\lambda }\right) ^{j}\frac{d^{\frac{j}{2}}}{dt^{\frac{j}{2}}}\psi
_{k,\alpha }^{el}=\frac{\alpha }{\sqrt{2}}\left( 1-\psi _{k,\alpha
}^{el}\right) -\frac{c_{k}}{\sqrt{\pi t}},  \label{gam3}
\end{equation}%
\textit{where }$c_{k}=1$\textit{\ for }$k$\textit{\ odd\ and }$c_{k}=0$%
\textit{\ for }$k$\textit{\ even. The initial conditions are }$\psi
_{k,\alpha }^{el}(0)=1,$\textit{\ for any }$k\geq 1$\textit{\ and }%
\begin{eqnarray}
\left. \frac{d^{r}}{dt^{r}}\psi _{k,\alpha }^{el}(t)\right\vert _{t=0}
&=&0,\quad r=1,...,\frac{k-1}{2},\text{ \textit{for odd} }k>1  \label{gam4}
\\
\left. \frac{d^{r}}{dt^{r}}\psi _{k,\alpha }^{el}(t)\right\vert _{t=0}
&=&0,\quad r=1,...,\frac{k}{2}-1,\text{ \textit{for even} }k>1.  \notag
\end{eqnarray}

\noindent \textbf{Proof \ }We rewrite (\ref{gam}) as follows:%
\begin{equation*}
\mathcal{L}\left\{ \psi _{k,\alpha }^{el};\eta \right\} \eta (\sqrt{2\eta }%
+\alpha )\sum_{j=0}^{k}\binom{k}{j}2^{\frac{j}{2}}\lambda ^{k-j}\eta ^{\frac{%
j}{2}}=(\sqrt{2\eta }+\alpha )\sum_{j=0}^{k}\binom{k}{j}2^{\frac{j}{2}%
}\lambda ^{k-j}\eta ^{\frac{j}{2}}-\sqrt{2\eta }\lambda ^{k}
\end{equation*}%
so that we get%
\begin{eqnarray}
&&\sum_{j=0}^{k}\binom{k}{j}\left( \frac{\sqrt{2}}{\lambda }\right) ^{j}%
\left[ \widetilde{\psi }_{k,\alpha }^{el}\eta ^{\frac{j}{2}+\frac{1}{2}%
}-\eta ^{\frac{j}{2}-\frac{1}{2}}\right] +\frac{\alpha }{\sqrt{2}}%
\sum_{j=1}^{k}\binom{k}{j}\left( \frac{\sqrt{2}}{\lambda }\right) ^{j}\left[
\widetilde{\psi }_{k,\alpha }^{el}\eta ^{\frac{j}{2}}-\eta ^{\frac{j}{2}-1}%
\right]  \notag \\
&=&\frac{\alpha }{\sqrt{2}}\left[ \frac{1}{\eta }-\widetilde{\psi }%
_{k,\alpha }^{el}\right] -\frac{1}{\sqrt{\eta }},  \label{gam5}
\end{eqnarray}%
where we have denoted $\widetilde{\psi }_{k,\alpha }^{el}=\mathcal{L}\left\{
\psi _{k,\alpha }^{el};\eta \right\} $ for brevity. From the Laplace
transform (\ref{gam5}), by taking into account (\ref{jen}) and the initial
conditions (\ref{gam4}), we can obtain equation (\ref{gam3}) with $c_{k}=1$.
For the initial conditions (\ref{gam4}) we use an argument similar to that
of Theorem 2.4, with the only additional care that, in the case of even $k$,
the highest order derivative, i.e. $\frac{d^{\frac{k}{2}}}{dt^{\frac{k}{2}}}%
\psi _{k,\alpha }^{el}$, does not vanish in $t=0$, as can be ascertained by
applying (\ref{derr}) to (\ref{gam1}): indeed we get%
\begin{equation*}
\left. \frac{d^{\frac{k}{2}}}{dt^{\frac{k}{2}}}\psi _{k,\alpha
}^{el}(t)\right\vert _{t=0}=\left. -\left( \frac{\lambda }{\sqrt{2}}\right)
^{k}\sum_{l=0}^{\infty }\left( -\frac{\alpha }{\sqrt{2}}\right) ^{l}t^{\frac{%
l}{2}}E_{\frac{1}{2},\frac{l}{2}+1}^{k}\left( -\frac{\lambda \sqrt{t}}{\sqrt{%
2}}\right) \right\vert _{t=0}=-\left( \frac{\lambda }{\sqrt{2}}\right) ^{k}.
\end{equation*}%
Therefore formula (\ref{gam5}), for even $k$, must be modified as follows%
\begin{eqnarray*}
&&\sum_{j=0}^{k-1}\binom{k}{j}\left( \frac{\sqrt{2}}{\lambda }\right) ^{j}%
\left[ \widetilde{\psi }_{k,\alpha }^{el}\eta ^{\frac{j}{2}+\frac{1}{2}%
}-\eta ^{\frac{j}{2}-\frac{1}{2}}\right] +\left( \frac{\sqrt{2}}{\lambda }%
\right) ^{k}\left[ \widetilde{\psi }_{k,\alpha }^{el}\eta ^{\frac{k}{2}+%
\frac{1}{2}}-\eta ^{\frac{k}{2}-\frac{1}{2}}-\frac{1}{\sqrt{\eta }}\left.
\frac{d^{\frac{k}{2}}}{dt^{\frac{k}{2}}}\psi _{k,\alpha }^{el}(t)\right\vert
_{t=0}\right] + \\
&&+\frac{\alpha }{\sqrt{2}}\sum_{j=1}^{k}\binom{k}{j}\left( \frac{\sqrt{2}}{%
\lambda }\right) ^{j}\left[ \widetilde{\psi }_{k,\alpha }^{el}\eta ^{\frac{j%
}{2}}-\eta ^{\frac{j}{2}-1}\right] \\
&=&\frac{\alpha }{\sqrt{2}}\left[ \frac{1}{\eta }-\widetilde{\psi }%
_{k,\alpha }^{el}\right] ,
\end{eqnarray*}%
so that we get (\ref{gam3}), with $c_{k}=0.$

As a further check, it is easy to see that, for $k=1,$ the latter reduces to
equation (\ref{el9}).\hfill $\blacksquare $

\section{Fractional relaxation equation of distributed order}

We consider now an extension of the fractional relaxation equation (\ref{due}%
) obtained by adding the hypothesis that the fractional order $\nu $ is not
a constant but a random variable with distribution $n(\nu )$. Thus we will
study the \textit{distributed order fractional relaxation equation} defined
as%
\begin{equation}
\int_{0}^{1}\frac{d^{\nu }\psi }{dt^{\nu }}n(\nu )d\nu =-\lambda \psi ,\quad
t>0,  \label{uno.2}
\end{equation}%
where, by assumption,
\begin{equation}
n(\nu )\geq 0\text{,\qquad }\int_{0}^{1}n(\nu )d\nu =1,\quad \nu \in \left(
0,1\right] ,  \label{bi1}
\end{equation}%
subject to the initial condition $\psi (0)=1$. As a special case, for $n(\nu
)=\delta (\nu -\overline{\nu })$ and a particular value of $\overline{\nu }%
\in \left( 0,1\right) ,$ equation (\ref{uno.2}) reduces to (\ref{due}).

We adopt here the following particular form for the density of the
fractional order $\nu $:%
\begin{equation}
n(\nu )=n_{1}\delta (\nu -\nu _{1})+n_{2}\delta (\nu -\nu _{2}),\qquad 0<\nu
_{1}<\nu _{2}\leq 1,  \label{bi3}
\end{equation}%
for $n_{1},n_{2}\geq 0$ and such that $n_{1}+n_{2}=1$ (conditions (\ref{bi1}%
) are trivially fulfilled). The density (\ref{bi3}) has been already used by
\cite{mapagn} and \cite{che1}, in the analysis of the so-called double-order
time-fractional diffusion equation, and corresponds to the case of a
subdiffusion with retardation. Moreover, it was applied in \cite{beg3} in
the context of recursive equations of fractional order, where the equation
governing the Poisson process has been extended by introducing two
fractional time derivatives.

Under assumption (\ref{bi3}), equation (\ref{uno.2}) becomes%
\begin{equation}
n_{1}\frac{d^{\nu _{1}}}{dt^{\nu _{1}}}\psi +n_{2}\frac{d^{\nu _{2}}}{%
dt^{\nu _{2}}}\psi =-\lambda \psi ,\quad t>0  \label{bi2}
\end{equation}%
and the corresponding solution $\psi _{\nu _{1},\nu _{2}}$ coincides with
the so-called \textit{double-order fractional relaxation} studied by \cite%
{mai}. They provide for $\psi _{\nu _{1},\nu _{2}}$ an integral expression
and some asymptotic representations. We present here an analytic form of the
fundamental solution to (\ref{bi2}) in terms of GML\ functions as well as a
probabilistic representation in terms of crossing probabilities, in line
with the results of the previous sections.

\

\noindent \textbf{Theorem 3.1 }\textit{The solution to equation (\ref{bi2})
with the initial condition }$\psi (0)=1$\textit{\ can be written as follows:}%
\begin{equation}
\psi _{\nu _{1},\nu _{2}}(t)=1-\frac{\lambda t^{\nu _{2}}}{n_{2}}%
\sum_{r=0}^{\infty }\left( -\frac{n_{1}t^{\nu _{2}-\nu _{1}}}{n_{2}}\right)
^{r}E_{\nu _{2},\nu _{2}+(\nu _{2}-\nu _{1})r+1}^{r+1}\left( -\frac{\lambda
t^{\nu _{2}}}{n_{2}}\right) .  \label{bi4}
\end{equation}%
\textit{\ }

\noindent \textbf{Proof }By taking the Laplace transform of (\ref{bi2}) we
get%
\begin{equation}
n_{1}\eta ^{\nu _{1}}\mathcal{L}\left\{ \psi _{\nu _{1},\nu _{2}};\eta
\right\} -\eta ^{\nu _{1}}+n_{2}\eta ^{\nu _{2}}\mathcal{L}\left\{ \psi
_{\nu _{1},\nu _{2}};\eta \right\} -\eta ^{\nu _{2}}=-\lambda \mathcal{L}%
\left\{ \psi _{\nu _{1},\nu _{2}};\eta \right\} ,  \label{bi5}
\end{equation}%
whose solution can be written as%
\begin{eqnarray*}
\mathcal{L}\left\{ \psi _{\nu _{1},\nu _{2}};\eta \right\} &=&\frac{%
n_{1}\eta ^{\nu _{1}}+n_{2}\eta ^{\nu _{2}}}{\eta (\lambda +n_{1}\eta ^{\nu
_{1}}+n_{2}\eta ^{\nu _{2}})} \\
&=&\frac{1}{\eta }-\frac{\lambda }{\eta }\frac{1}{\lambda +n_{2}\eta ^{\nu
_{2}}}\frac{1}{1+\frac{n_{1}\eta ^{\nu _{1}}}{\lambda +n_{2}\eta ^{\nu _{2}}}%
} \\
&=&\frac{1}{\eta }-\frac{\lambda }{\eta }\frac{1}{\lambda +n_{2}\eta ^{\nu
_{2}}}\sum_{r=0}^{\infty }\left( -\frac{n_{1}\eta ^{\nu _{1}}}{\lambda
+n_{2}\eta ^{\nu _{2}}}\right) ^{r} \\
&=&\frac{1}{\eta }-\frac{\lambda }{n_{2}}\sum_{r=0}^{\infty }\left( -\frac{%
n_{1}}{n_{2}}\right) ^{r}\frac{\eta ^{\nu _{1}r-1}}{\left( \eta ^{\nu _{2}}+%
\frac{\lambda }{n_{2}}\right) ^{r+1}}.
\end{eqnarray*}%
By applying formula (\ref{pra}), we easily get (\ref{bi4}). As a check we
can see that (\ref{bi4}) reduces to (\ref{pr2}), for $n_{1}=0,$ $n_{2}=1$, $%
\nu _{2}=\nu $, since equation (\ref{bi2}) becomes, in this case, the
fractional relaxation equation (\ref{due}).\hfill $\blacksquare $

\

Despite the apparent similarity of (\ref{bi4}) with (\ref{gam1}), they are
deeply different: while for $\psi _{\nu _{1},\nu _{2}}$ the sum is extended
to the third (upper) parameter of the GML function, this is not the case for
$\psi _{k,\alpha }^{el}.$ This is also reflected in the asymptotic behavior
of the fractional relaxation of distributed order, which does not deviate
from the usual relaxation behavior (unlike $\psi _{k,\alpha }^{el}$). We can
study the limit directly from (\ref{bi4}), by applying formula (\ref{asy}),
as follows%
\begin{eqnarray}
&&\psi _{\nu _{1},\nu _{2}}(t) \\
&=&1-\frac{\lambda }{n_{2}}\sum_{r=0}^{\infty }\left( -\frac{n_{1}}{n_{2}}%
\right) ^{r}\frac{1}{2\pi i}\int_{0}^{\infty }e^{-zt}z^{\nu _{1}r-1}\left[
\frac{e^{-i\pi \nu _{2}-i\pi (\nu _{2}-\nu _{1})r}}{\left( z^{\nu _{2}}+%
\frac{\lambda }{n_{2}}e^{-i\pi \nu _{2}}\right) ^{r+1}}-\frac{e^{i\pi \nu
_{2}+i\pi (\nu _{2}-\nu _{1})r}}{\left( z^{\nu _{2}}+\frac{\lambda }{n_{2}}%
e^{i\pi \nu _{2}}\right) ^{r+1}}\right] .  \notag
\end{eqnarray}%
Thus, for $t\rightarrow 0,$ we get
\begin{eqnarray}
&&\psi _{\nu _{1},\nu _{2}}(t)=1-\frac{\lambda }{n_{2}}\sum_{r=0}^{\infty
}\left( -\frac{n_{1}t^{\nu _{2}-\nu _{1}}}{n_{2}}\right) ^{r}\frac{t^{\nu
_{2}}}{2\pi i}\int_{0}^{\infty }e^{-w}w^{\nu _{1}r-1}\cdot  \label{dop} \\
&&\cdot \left[ \frac{e^{-i\pi \nu _{2}-i\pi (\nu _{2}-\nu _{1})r}}{\left(
w^{\nu _{2}}+\frac{\lambda t^{\nu _{2}}}{n_{2}}e^{-i\pi \nu _{2}}\right)
^{r+1}}-\frac{e^{i\pi \nu _{2}+i\pi (\nu _{2}-\nu _{1})r}}{\left( w+\frac{%
\lambda t^{\nu _{2}}}{n_{2}}e^{i\pi \nu _{2}}\right) ^{r+1}}\right]  \notag
\\
&\simeq &1-\frac{\lambda t^{\nu _{2}}}{n_{2}}\sum_{r=0}^{\infty }\left( -%
\frac{n_{1}t^{\nu _{2}-\nu _{1}}}{n_{2}}\right) ^{r}\frac{\sin (-\pi (\nu
_{1}r-\nu _{2}r-\nu _{2}))}{\pi }\Gamma \left( \nu _{1}r-\nu _{2}r-\nu
_{2}\right)  \notag \\
&=&\left[ \text{by the reflection property of the Gamma function}\right]
\notag \\
&=&1-\frac{\lambda t^{\nu _{2}}}{n_{2}}\sum_{r=0}^{\infty }\left( -\frac{%
n_{1}t^{\nu _{2}-\nu _{1}}}{n_{2}}\right) ^{r}\frac{1}{\Gamma \left( 1+\nu
_{2}r+\nu _{2}-\nu _{1}r\right) }=1-\frac{\lambda t^{\nu _{2}}}{n_{2}}\frac{1%
}{\Gamma \left( 1+\nu _{2}\right) }+o(t^{\nu _{2}}),  \notag
\end{eqnarray}%
while, for $t\rightarrow \infty ,$ we analogously have that%
\begin{eqnarray}
&&\psi _{\nu _{1},\nu _{2}}(t)=1-\frac{\lambda }{n_{2}}\sum_{r=0}^{\infty
}\left( -\frac{n_{1}}{n_{2}t^{\nu _{1}}}\right) ^{r}\frac{1}{2\pi i}%
\int_{0}^{\infty }e^{-w}w^{\nu _{1}r-1}\cdot  \label{dop2} \\
&&\cdot \left[ \frac{e^{-i\pi \nu _{2}-i\pi (\nu _{2}-\nu _{1})r}}{\left(
\left( \frac{w}{t}\right) ^{\nu _{2}}+\frac{\lambda }{n_{2}}e^{-i\pi \nu
_{2}}\right) ^{r+1}}-\frac{e^{i\pi \nu _{2}+i\pi (\nu _{2}-\nu _{1})r}}{%
\left( \left( \frac{w}{t}\right) ^{\nu _{2}}+\frac{\lambda }{n_{2}}e^{i\pi
\nu _{2}}\right) ^{r+1}}\right]  \notag \\
&\simeq &1-\sum_{r=0}^{\infty }\left( -\frac{n_{1}}{\lambda t^{\nu _{1}}}%
\right) ^{r}\frac{\sin (\pi \nu _{1}r)}{\pi }\Gamma \left( \nu _{1}r\right)
\notag \\
&=&1-\sum_{r=0}^{\infty }\left( -\frac{n_{1}}{\lambda t^{\nu _{1}}}\right)
^{r}\frac{1}{\Gamma \left( 1-\nu _{1}r\right) }=\frac{n_{1}}{\lambda t^{\nu
_{1}}}\frac{1}{\Gamma \left( 1-\nu _{1}\right) }+o(t^{-\nu _{1}}).  \notag
\end{eqnarray}%
The previous expressions coincides with formula (4.16) of \cite{mai}, which
has been obtained in a different way, directly from the Laplace transform of
$\psi _{\nu _{1},\nu _{2}}.$

\

We present now a probabilistic form of the solution $\psi _{\nu _{1},\nu
_{2}}$, which is in line with the analysis carried out so far, in terms of
crossing probability of a random boundary by a stochastic process, that will
be denoted, in this case, by $\mathcal{T}_{\nu _{1},\nu _{2}}(t),t>0.$ To
this aim we will compare equation (\ref{bi2}) with the equation governing
the probabilities $\widetilde{p}_{k}$ of the distributed order fractional
Poisson process $\mathcal{N}_{\nu _{1},\nu _{2}}(t),t>0$ studied in \cite%
{beg3}, i.e.%
\begin{equation}
\int_{0}^{1}\frac{d^{\nu }p_{k}}{dt^{\nu }}n(\nu )d\nu =-\lambda
(p_{k}-p_{k-1}),\qquad k\geq 0,\text{ }p_{-1}(t)=0  \label{dop3}
\end{equation}%
Indeed (\ref{uno.2}) can be considered a special case of (\ref{dop3}) for $%
k=0$ and, if we add the assumption (\ref{bi3}), we get (\ref{bi2}). Thus we
can use the results proved in \cite{beg3} and write that%
\begin{equation}
\psi _{\nu _{1},\nu _{2}}(t)=\widetilde{p}_{0}(t)=\Pr \left\{ \mathcal{N}%
_{\nu _{1},\nu _{2}}(t)=0\right\} =\Pr \left\{ N(\mathcal{T}_{\nu _{1},\nu
_{2}}(t))=0\right\}  \label{dop5}
\end{equation}%
where $N$ is the standard Poisson process (with intensity $\lambda $) and $%
\mathcal{T}_{\nu _{1},\nu _{2}}$ is a random process (independent from $N$)
with density
\begin{equation}
q_{\nu _{1},\nu _{2}}(y,t)=n_{1}\int_{0}^{t}\overline{p}_{\nu
_{2}}(t-s;y)q_{\nu _{1}}(y,s)ds+n_{2}\int_{0}^{t}\overline{p}_{\nu
_{1}}(t-s;y)q_{\nu _{2}}(y,s)ds.  \label{dop4}
\end{equation}%
In (\ref{dop4}) $\overline{p}_{\nu _{j}}(\cdot ;z)$\ denotes the density of
a stable random variable $X_{\nu _{j}}$\ of index $\nu _{j}\in \left( 0,1%
\right] ,$\ for $j=1,2,$\ with parameters equal $\beta =1,$\ $\mu =0$\ and $%
\sigma =\left( n_{j}|y|\cos \frac{\pi \nu _{j}}{2}\right) ^{1/\nu _{j}}$ and
$q_{\nu _{j}}$, for $j=1,2,$\ was defined in (\ref{qu}). Another form of the
density $q_{\nu _{1},\nu _{2}}$ is given by the following series expression%
\begin{eqnarray}
&&q_{\nu _{1},\nu _{2}}(y,t)  \label{bi11} \\
&=&\frac{n_{1}}{\lambda t^{\nu _{1}}}\sum_{r=0}^{\infty }\frac{1}{r!}\left( -%
\frac{n_{2}|y|}{\lambda t^{\nu _{2}}}\right) ^{r}\mathcal{W}_{-\nu
_{1},1-\nu _{2}r-\nu _{1}}\left( -\frac{n_{1}|y|}{\lambda t^{\nu _{1}}}%
\right) +  \notag \\
&&+\frac{n_{2}}{\lambda t^{\nu _{2}}}\sum_{r=0}^{\infty }\frac{1}{r!}\left( -%
\frac{n_{1}|y|}{\lambda t^{\nu _{1}}}\right) ^{r}\mathcal{W}_{-\nu
_{2},1-\nu _{1}r-\nu _{2}}\left( -\frac{n_{2}|y|}{\lambda t^{\nu _{2}}}%
\right) .  \notag
\end{eqnarray}%
From (\ref{dop5}) we get
\begin{equation}
\psi _{\nu _{1},\nu _{2}}(t)=\int_{0}^{\infty }e^{-\lambda y}q_{\nu _{1},\nu
_{2}}(y,t)dy=\Pr \left\{ \mathcal{T}_{\nu _{1},\nu _{2}}(t)<U\right\} .
\label{dop6}
\end{equation}%
It is also proved in \cite{beg3} that the transition density $q_{\nu
_{1},\nu _{2}}$ coincides with the folded solution
\begin{equation}
q_{\nu _{1},\nu _{2}}(y,t)=\left\{
\begin{array}{l}
2v(y,t),\qquad y\geq 0 \\
0,\qquad y<0%
\end{array}%
\right.  \label{ops}
\end{equation}%
of the following fractional diffusion equation

\begin{equation}
\left( n_{1}\frac{\partial ^{\nu _{1}}v}{\partial t^{\nu _{1}}}+n_{2}\frac{%
\partial ^{\nu _{2}}v}{\partial t^{\nu _{2}}}\right) ^{2}=\frac{\partial
^{2}v}{\partial y^{2}},\quad y\in \mathbb{R},t>0,\;n_{1},n_{2}>0,  \label{e3}
\end{equation}%
for $0<\nu _{1}<\nu _{2}\leq 1$,\ with initial conditions%
\begin{equation}
\left\{
\begin{array}{l}
v(y,0)=\delta (y),\text{ for }0<\nu _{1}<\nu _{2}\leq 1 \\
\left. \frac{\partial }{\partial t}v(y,t)\right\vert _{t=0}=0\text{ for }%
\frac{1}{2}<\nu _{1}<\nu _{2}\leq 1%
\end{array}%
\right. .  \label{e5}
\end{equation}%
In alternative to (\ref{e3})-(\ref{e5}) it can be proved (as we will see
below in a special case) that $q_{\nu _{1},\nu _{2}}$ solves also the other
equation%
\begin{equation}
n_{1}\frac{\partial ^{\nu _{1}}v}{\partial t^{\nu _{1}}}+n_{2}\frac{\partial
^{\nu _{2}}v}{\partial t^{\nu _{2}}}=-\frac{\partial v}{\partial y},\quad
y,t>0,\;n_{1},n_{2}>0,\text{ }v(y,0)=\delta (y),  \label{bi12}
\end{equation}%
which is the distributed order analogue of (\ref{nofold}). In order to get a
more explicit expression of the density $q_{\nu _{1},\nu _{2}},$ we consider
the special, but relevant, case where $\nu _{1}=\frac{1}{2}$ and $\nu _{2}=1$%
.

\

\noindent \textbf{Theorem 3.2 }\textit{The solution to the fractional
relaxation equation }%
\begin{equation}
n_{1}\frac{d^{1/2}\psi }{dt^{1/2}}+n_{2}\frac{d\psi }{dt}=-\lambda \psi
,\quad t>0,  \label{dop7}
\end{equation}%
\textit{\ with the initial condition }$\psi (0)=1,$\textit{\ can be
expressed as follows:}%
\begin{equation}
\psi _{\frac{1}{2},1}(t)=\Pr \left\{ \mathcal{T}_{\frac{1}{2}%
,1}(t)<U\right\} ,  \label{dop12}
\end{equation}%
\textit{where }$U$\textit{\ is an exponential r.v. with parameter }$\lambda $%
\textit{\ and the transition density of }$\mathcal{T}_{\frac{1}{2}%
,1}(t),t>0, $\emph{\ }\textit{is given by}%
\begin{equation}
q_{\frac{1}{2},1}(y,t)=\frac{n_{1}(t-\frac{n_{2}}{2}y)}{\sqrt{\pi }}\frac{%
e^{-\frac{n_{1}^{2}y^{2}}{4(t-n_{2}y)}}}{\sqrt{(t-n_{2}y)^{3}}},\qquad t>0,%
\text{ }0<y<\frac{t}{n_{2}},  \label{dop8}
\end{equation}%
\textit{and satisfies the fractional equation}%
\begin{equation}
n_{1}\frac{\partial ^{1/2}q}{\partial t^{1/2}}+n_{2}\frac{\partial q}{%
\partial t}=-\frac{\partial q}{\partial y},\qquad q(y,0)=\delta (y).
\label{dop10}
\end{equation}%
\textbf{Proof }\ It has been proved in \cite{beg3} that for $\nu _{2}=1$ and
$\nu _{1}=\nu \in (0,1)$ the density (\ref{dop4}), can be expressed as%
\begin{equation}
q_{\nu ,1}(y,t)=n_{1}I^{\nu }(\overline{\overline{p}}_{\nu }(\cdot
;y))(t)+n_{2}\overline{\overline{p}}_{\nu }(t;y),  \label{dop9}
\end{equation}%
where $I^{\nu }$ is the Riemann-Liouville fractional integral of order $\nu $
and $\overline{\overline{p}}_{\nu }$ denotes a stable law of index $\nu $
and parameters equal to $\beta =1,$ $\mu =n_{2}|y|,$ $\sigma =\left(
n_{1}|y|\cos \frac{\pi \nu }{2}\right) ^{1/\nu }.$ If we put moreover $\nu
=1/2,$ we can recognize in $\overline{\overline{p}}_{\frac{1}{2}}$ the Lévy
distribution, so that the density (\ref{dop9}) becomes
\begin{eqnarray*}
&&q_{\frac{1}{2},1}(y,t) \\
&=&\frac{n_{1}}{\sqrt{\pi }}\int_{0}^{t}(t-s)^{-\frac{1}{2}}\overline{%
\overline{p}}_{\frac{1}{2}}(s;y)ds+n_{2}\overline{\overline{p}}_{\frac{1}{2}%
}(t;y) \\
&=&\frac{n_{1}^{2}y}{2\pi }\int_{n_{2}y}^{t}(t-s)^{-\frac{1}{2}}\frac{e^{-%
\frac{n_{1}^{2}y^{2}}{4(s-n_{2}y)}}}{\sqrt{(s-n_{2}y)^{3}}}ds+\frac{%
n_{1}n_{2}y}{2\pi }\frac{e^{-\frac{n_{1}^{2}y^{2}}{4(t-n_{2}y)}}}{\sqrt{%
(t-n_{2}y)^{3}}}\mathit{1}_{\left\{ 0<y<\frac{t}{n_{2}}\right\} } \\
&=&\frac{n_{1}^{2}y}{2\pi }\int_{0}^{t-n_{2}y}(t-n_{2}y-z)^{-\frac{1}{2}}%
\frac{e^{-\frac{n_{1}y^{2}}{4z}}}{\sqrt{z^{3}}}dz+\frac{n_{1}n_{2}y}{2\pi }%
\frac{e^{-\frac{n_{1}^{2}y^{2}}{4(t-n_{2}y)}}}{\sqrt{(t-n_{2}y)^{3}}}\mathit{%
1}_{\left\{ 0<y<\frac{t}{n_{2}}\right\} } \\
&=&\left[ \text{by the identity (3.8) of \cite{ors}}\right] \\
&=&\left[ \frac{n_{1}e^{-\frac{n_{1}^{2}y^{2}}{4(s-n_{2}y)}}}{\sqrt{\pi
(t-n_{2}y)}}+\frac{n_{1}n_{2}y}{2\pi }\frac{e^{-\frac{n_{1}^{2}y^{2}}{%
4(t-n_{2}y)}}}{\sqrt{(t-n_{2}y)^{3}}}\right] \mathit{1}_{\left\{ 0<y<\frac{t%
}{n_{2}}\right\} },
\end{eqnarray*}%
which coincides with (\ref{dop8}). In order to show that the latter
satisfies the fractional relaxation equation (\ref{dop10}), we evaluate its
Laplace transform, which reads:%
\begin{eqnarray}
&&\mathcal{L}\left\{ q_{\frac{1}{2},1}(y,\cdot );\eta \right\}  \label{dop11}
\\
&=&\int_{n_{2}y}^{\infty }\frac{n_{1}(t-\frac{n_{2}}{2}y)}{\sqrt{\pi }}\frac{%
e^{-\frac{n_{1}^{2}y^{2}}{4(t-n_{2}y)}-\eta t}}{\sqrt{(t-n_{2}y)^{3}}}dt
\notag \\
&=&-\frac{n_{1}}{\sqrt{\pi }}\frac{\partial }{\partial \eta }\left\{
\int_{n_{2}y}^{\infty }\frac{e^{-\frac{n_{1}^{2}y^{2}}{4(t-n_{2}y)}-\eta t}}{%
\sqrt{(t-n_{2}y)^{3}}}dt\right\} -\frac{n_{1}n_{2}y}{2\sqrt{\pi }}%
\int_{0}^{\infty }e^{-\eta z-\eta n_{2}y}\frac{e^{-\frac{n_{1}^{2}y^{2}}{4z}}%
}{\sqrt{z^{3}}}dz  \notag \\
&=&-\frac{\partial }{\partial \eta }\left\{ \frac{2}{y}e^{-\eta n_{2}y-\sqrt{%
\eta }n_{1}y}dt\right\} -n_{2}e^{-\eta n_{2}y-\sqrt{\eta }n_{1}y}  \notag \\
&=&\left( n_{2}+n_{1}\eta ^{-1/2}\right) e^{-(n_{2}\eta +n_{1}\eta ^{1/2})y}.
\notag
\end{eqnarray}%
In (\ref{dop11}) we have applied the well-known formula of the Laplace
transform of the first-passage time of a Brownian motion. It is easy to
check that
\begin{equation*}
\int_{0}^{\infty }e^{-\lambda y}\mathcal{L}\left\{ q_{\frac{1}{2},1}(y,\cdot
);\eta \right\} dy=\frac{n_{2}+n_{1}\eta ^{-1/2}}{n_{2}\eta +n_{1}\eta
^{1/2}+\lambda },
\end{equation*}%
which is equal to the Laplace transform of $\psi _{\nu _{1},\nu _{2}},$ for $%
\nu _{1}=1/2$ and $\nu _{2}=1$ (given in Theorem 2.6 of \cite{beg3}), thus
proving result (\ref{dop12}). If we now take the Fourier transform of (\ref%
{dop11}) we get%
\begin{eqnarray}
&&\mathcal{F}\left\{ \mathcal{L}\left\{ q_{\frac{1}{2},1};\eta \right\}
;\beta \right\} =\int_{0}^{\infty }e^{i\beta y}\mathcal{L}\left\{ q_{\frac{1%
}{2},1}(y,\cdot );\eta \right\} dy  \label{dop15} \\
&=&\left( n_{2}+n_{1}\eta ^{-1/2}\right) \int_{0}^{\infty }e^{i\beta
y}e^{-(n_{2}\eta +n_{1}\eta ^{1/2})y}dy  \notag \\
&=&\frac{n_{2}+n_{1}\eta ^{-1/2}}{n_{2}\eta +n_{1}\eta ^{1/2}+i\beta },
\notag
\end{eqnarray}%
which coincides with the solution to equation (\ref{dop10}) converted, via
Laplace-Fourier transform, into%
\begin{equation*}
(n_{1}\eta ^{1/2}+n_{2}\eta )\mathcal{L}\left\{ q_{\frac{1}{2},1}(y,\cdot
);\eta \right\} -(n_{1}\eta ^{-1/2}+n_{2})\delta (y)=-\frac{\partial }{%
\partial y}\mathcal{L}\left\{ q_{\frac{1}{2},1}(y,\cdot );\eta \right\}
\end{equation*}%
and%
\begin{equation*}
(n_{1}\eta ^{1/2}+n_{2}\eta +i\beta )\mathcal{F}\left\{ \mathcal{L}\left\{
q_{\frac{1}{2},1};\eta \right\} ;\beta \right\} =(n_{1}\eta ^{-1/2}+n_{2}).
\end{equation*}%
From (\ref{dop15}) it is evident that (\ref{dop8}) is well-defined and
integrates to one, since for $\beta =0$ we get $1/\eta .$\hfill $%
\blacksquare $

\

\noindent \textbf{Remark 3.1 }If we consider the two opposite special cases $%
n_{2}=0$ and $n_{1}=0,$ the trajectories of the process $\mathcal{T}_{\frac{1%
}{2},1}$ can be considered as \textquotedblleft
interpolation\textquotedblright\ between those of a free reflecting Brownian
motion and the straight line $y=t/n_{2}.$ Indeed in the first case the
density (\ref{dop8}) becomes%
\begin{equation*}
q_{\frac{1}{2},1}(y,t)=\frac{n_{1}e^{-\frac{n_{1}^{2}y^{2}}{4t}}}{\sqrt{\pi t%
}},\qquad y,t>0,
\end{equation*}%
while in the second we can write (\ref{dop9}) as $q_{\frac{1}{2}%
,1}(y,t)=n_{2}\overline{\overline{p}}_{\frac{1}{2}}(t;y)=n_{2}\delta
(t-n_{2}y)$, since in this case $\sigma =0.$ It is evident from (\ref{dop8})
that the trajectories of $\mathcal{T}_{\frac{1}{2},1}$, for any $n_{1}$,$%
n_{2}>0$ are forced under the line $y=t/n_{2}$ and this is reflected in the
asymptotic behavior of the crossing probability $\psi _{\frac{1}{2},1}$,
which can be deduced from (\ref{dop}) and (\ref{dop2}) and summed up as
follows:%
\begin{equation}
\psi _{\frac{1}{2},1}(t)\simeq \left\{
\begin{array}{l}
1-\frac{\lambda t}{n_{2}},\qquad 0<t<<1 \\
\frac{n_{1}}{\lambda \sqrt{\pi t}},\qquad t\rightarrow \infty%
\end{array}%
\right. .  \label{dop13}
\end{equation}%
By comparing (\ref{dop13}) with (\ref{asy4}) we can conclude that $\psi _{%
\frac{1}{2},1}$ displays the same limiting behavior of $\psi _{\frac{1}{2}%
}(t)=\Pr \left\{ |B(t)|<U\right\} $, for $t\rightarrow \infty .$ On the
contrary, for $t\rightarrow 0$, it behaves as the standard relaxation (up to
a constant) and thus tends to one much faster than $\psi _{\frac{1}{2},1}.$
We recall that similar limiting features were exhibited by the crossing
probability $\psi ^{+}$ of the Brownian sojourn time process (see (\ref{psi3}%
)).

\

For the reader's convenience we summarize the limiting behavior of the
crossing probabilities analyzed in the previous sections in the following
tables:

\

\begin{center}
\textbf{Table 1}: Limiting behavior for $t\rightarrow 0$%
\begin{equation*}
\begin{array}{c}
\psi (t)\simeq 1-\lambda t \\
\psi _{\nu }(t)\simeq 1-\frac{\lambda t^{\nu }}{\Gamma (1+\nu )} \\
\psi ^{+}(t)\simeq 1-\frac{\lambda t}{2} \\
\psi ^{T}(t)\simeq 1-\sqrt{2\lambda }t \\
\psi ^{\gamma }(t)\simeq \frac{1}{(1+2\lambda t)^{\gamma /2}} \\
\psi ^{el}(t)\simeq 1-\frac{\lambda \sqrt{2t}}{\sqrt{\pi }} \\
\psi _{\frac{1}{2}}^{k}(t)\simeq 1-\frac{(\lambda \sqrt{t})^{k}}{\Gamma
\left( \frac{k}{2}+1\right) } \\
\psi _{k,\alpha }^{el}(t)\simeq 1-\frac{(\frac{\lambda \sqrt{t}}{2})^{k}}{%
\Gamma \left( \frac{k}{2}+1\right) } \\
\psi _{\frac{1}{2},1}(t)\simeq 1-\frac{\lambda t}{n_{2}}%
\end{array}%
\end{equation*}
\end{center}

\

\begin{center}
\textbf{Table 2}: Limiting behavior for $t\rightarrow \infty $%
\begin{equation*}
\begin{array}{c}
\psi (t)\simeq e^{-\lambda t} \\
\psi _{\nu }(t)\simeq \frac{1}{\lambda t^{\nu }\Gamma (1-\nu )} \\
\psi ^{+}(t)\simeq \frac{1}{\sqrt{\lambda \pi t}} \\
\psi ^{T}(t)\simeq e^{-\sqrt{2\lambda }t} \\
\psi ^{\gamma }(t)\simeq \frac{1}{(1+2\lambda t)^{\gamma /2}} \\
\psi ^{el}(t)\simeq 1-\frac{\sqrt{2}}{\alpha \sqrt{\pi t}} \\
\psi _{\frac{1}{2}}^{k}(t)\simeq \frac{k}{\lambda \sqrt{\pi t}} \\
\psi _{k,\alpha }^{el}(t)\simeq 1-\frac{\sqrt{2}}{\alpha \sqrt{\pi t}} \\
\psi _{\frac{1}{2},1}(t)\simeq \frac{n_{1}}{\lambda \sqrt{\pi t}}%
\end{array}%
\end{equation*}
\end{center}

\

\end{document}